

\input amssym.def

\newsymbol\restrictedto 1316

\input epsf.sty


\font\bigbigrm = cmr10 scaled \magstep2

\def\AA{{\cal A}}
\def\CC{{\cal C}}
\def\EE{{\cal E}}
\def\FF{{\cal F}}
\def\GG{{\cal G}}
\def\UU{{\cal U}}
\def\Union{\bigcup}
\def\union{\cup}
\def\intersect{\cap}
\def\includedin{\subseteq}
\def\including{\supseteq}
\def\includes{\supseteq}
\def\proof#1:{\par\smallskip\noindent {\sl Proof #1\/}:\par\nobreak}
\def\ceiling#1{\lceil#1\rceil}
\def\iso{\simeq}
\let\<\langle
\let\>\rangle
\let\implies\Rightarrow

\def\section#1.#2\par{\par\medskip\noindent{\bf \S#1}.%
\quad \sl#2\par\smallskip\rm}
\def\introduction{\par\smallskip\noindent Introduction\par\smallskip}
\def\part#1.#2\par{\par\bigskip\noindent \bf #1.\quad #2\par\smallskip\rm}
\def\proclaim#1#2{\smallskip\par \noindent{\bf #1}\quad \sl #2}
\def\endproclaim{\par\smallskip \rm}
\def\qed{\hfill{\vrule height 8pt depth .05pt width 5pt}}

\def\vbar{{\bf v}}
\def\acl{\hbox{\rm acl}}
\def\cl{\hbox{\rm cl}}

\centerline{\bigbigrm \bf UNIVERSAL GRAPHS WITH FORBIDDEN SUBGRAPHS AND
ALGEBRAIC CLOSURE}
\vskip 0.4in
\centerline{Gregory Cherlin
\footnote{${}^{1}$}{\it Research supported in part by NSF grant DMS
9501176}}
\centerline{Dept. Mathematics, Rutgers University, Busch Campus, New 
Brunswick, NJ 08903}
\vskip 0.1in
\centerline{Saharon Shelah
\footnote{${}^{2}$}{\it Research supported in part by U.S.-Israel
Binational Science Foundation grant 0377215}}
\centerline{Dept. Mathematics, Hebrew University, Jerusalem, Israel,
and Rutgers University, New Brunswick, NJ 08903}
\vskip 0.1in
\centerline{Niandong Shi}
\centerline{Dept. Mathematics, East Stroudsburg University, East 
Stroudsburg, PA 18301}
\vskip 0.4in
\centerline{ABSTRACT}

{\leftskip=30pt \rightskip=30pt 

We apply model theoretic methods to the problem of existence of countable
universal graphs with finitely many forbidden connected subgraphs.  We
show that to a large extent the question reduces to one of local
finiteness of an associated ``algebraic closure" operator (Theorem 3, \S3). The
main applications are new examples of universal graphs with forbidden
subgraphs (\S\S 7, 8, and 10) 
and simplified treatments of some previously known
cases (\S\S 6.2,6.3). \par}

\baselineskip = 16pt

\introduction

We are concerned here with the following problem: give a finite set 
$\CC$ of finite connected graphs, does the class $\GG_{\CC}$ of
countable graphs which omit $\CC$ contain a universal element (one in
which all others are embeddable as induced subgraphs)?  Here we say that a
graph $G$ {\sl omits} a class $\CC$ of graphs if no graph in $\CC$
embeds as a subgraph of $G$.  The problem is to characterize those
classes
$\CC$ for which there is such a universal graph.  A more fundamental
problem is whether there is any effective characterization of these
classes
$\CC$; in other words, is there an algorithm which will produce the
answer in each case?  This problem remains open even when $\CC$
consists of a {\sl single} forbidden subgraph, though an accumulation of
evidence, some given in the present paper, suggests that at least this
instance should have an explicit and fairly simple solution.  We discuss
this further below.

Rado observed [Ra] that there is a universal countable graph.  This
corresponds to the case ${\CC}=\emptyset$.  Many other cases have been
considered in the literature [ChK,CS1,CS2,CST,FK1,Ko,FK2,GK,KMP,Ko,KP1,
KP2,Pa].
In particular [FK1] gives a complete solution for the case in which 
$\CC$ consists of a single 2-connected constraint, and [CST] treats the
case in which $\CC$ consists of a single tree with no vertex of degree
2.  

To date very few cases have been identified in which a universal countable
$\CC$-free graph exists.  For ${\CC}=\{C\}$ consisting of a single
constraint, the following cases are known to allow a universal graph: 
$C$ is complete; $C$ is a tree consisting of one path
to which at most one additional edge is attached; or $C$ is a ``bow-tie", a
particular graph of order 5.  We will add some additional families of
examples using model theoretic methods.

Another family of universal $\CC$-free graphs corresponds to the class
$\CC$ of odd cycles of length up to some specified bound
$(C_{2n+1}$-free graphs for $n\le N)$.  We generalize this to the case in
which $\CC$ is closed under homomorphism in an appropriate sense
(Theorem 4, and $\S 7$). In earlier work the positive results have
generally
come either from structure theorems for $\CC$-free graphs (notably in
the path-free case [KMP]) or from 
Fra\"{\i}ss\'e's amalgamation method, whose
complexity increases rapidly as the constraint class $\CC$ becomes more
complicated.  
Using the model theoretic notions of existential completeness and
algebraic closure for $\CC$-free graphs, we give a criterion for the
existence of a universal $\CC$-free graph which effectively short
circuits this process.  Our arguments can be converted into amalgamation
arguments in principle, but not in any very explicit way.  

Kojman conjectured in conversation years ago that closure of the constraint
class under homomorphic image might be a key condition in connection with
the existence of universal graphs. Our work confirms this conjecture in
one direction (Theorems 4 and 5)
and relates the condition directly to the broader issue
of the structure of the algebraic closure operator.

Our model theoretic methods are very close to those which have been used
in practice in all cases in which {\sl nonexistence} of universal graphs
has been established.  There they are typically referred to as ``rigidity"
arguments.  This amounts to a rephrasing in purely graph theoretic terms
of a more general model theoretic notion.  
We did this ourselves in [CST],
though in fact an awareness of the model theoretic  framework lay in the
background of the proof given there.  In our present work, we have reached 
the point at which such 
a reformulation of our methods would be counterproductive, as we
make use of general considerations which are well known in model theory
but have not yet played in explicit role in graph theory.  Accordingly,
the first half of the present paper lays the foundation of our approach,
recalling what we need from model theory and applying it in the case of
$\CC$-free graphs.  As we will see in Theorem 3 of $\S 3$, these ideas
produce much clearer results in the $\CC$-free context than one would
get in a more general model theoretic context.  This is really the key to
our whole analysis.

Applications of these general ideas are found in $\S \S 5-10$.  Most of the
cases considered in $\S\S 5,6,9$ were treated successfully in the past,
and are reexamined from our present point of view partly by way of
illustration and partly because our present viewpoint suggests
quantitative issues extending the earlier purely qualitative analysis.
That is, in cases in which our ``algebraic closure" operator is locally
finite, we consider its rate of growth.

New examples are given in $\S\S 7,8$.  In particular $\S 8$ is devoted to
an infinite family extending the ``bow-tie" example considered by 
Komj\'ath [Ko] using methods that have further potential. 
This is the hardest case treated here.

To conclude this introduction we take note of two directions which are
particularly promising for further work: the general problem of
effectivity, and the case of a single constraint.

\proclaim{Effectivity}
Given a finite set $\CC$ of finite connected graphs, determine
whether there is a universal countable $\CC$-free graph.
\endproclaim

It is by no means clear that this problem should have an effective
solution.  It is natural to consider a further generalization in which
graphs are replaced by vertex-colored, edge-colored, and directed graphs,
or more generally by relational structures for any finite relational
language . It seems likely however that this more general problem can
be reduced to the special case of graphs by a suitable encoding.  This is
one reason why the existence of an effective solution is doubtful, but at
the present time the question is entirely open.  

\proclaim{Single constraints}
Let $\UU$ be the collection of all finite connected graphs 
$C$ for which
there is a countable universal $C$-free graph, and let 
$$\UU_0=\{C\in \UU: \hbox{ every induced subgraph of $C$
 is in $\UU$}\}.$$
\endproclaim

\proclaim{Conjecture}
 $\UU = \UU_0$.
\endproclaim

While it is not likely that this conjecture will be proved {\sl a priori},
it may well turn out to be the case.  The point of the conjecture is that
it should be possible to determine $\UU_0$ explicitly using known
methods, and then rephrase the conjecture more explicitly.  In  
$\S 8$ below we will give a new family of examples in $\UU_0$.
What is
needed, apparently, is to continue that analysis, which will involve
substantial computation, and also to prove a number of further results on
nonexistence of universal graphs.  Note that by [FK1] any block (maximal
2-connected subgraph) of a graph in $\UU_0$ is complete, and the
results of [GK] can be combined with some similar constructions to reduce
the class of candidates for members of $\UU_0$ to graphs much like
those considered in $\S 8$.

\part I. General theory.

In this part we will discuss the application of some model
theoretic ideas to the general problem of the existence of universal
countable graphs with forbidden subgraphs.  In $\S 1$ we associate with a
class $\CC$ of finite graphs (usually taken to be connected) the class
$\GG_{\CC}$ of countable graphs ``omitting" $\CC$ and the
subclass ${\EE}_{\CC}$ of ``existentially complete" graphs in 
$\GG_{\CC}$. The key to the model theoretic approach is to
understand ${\EE}_{\CC}$.  In fact where a universal $\CC$-free
graph exists, it is often the case that ${\EE}_{\CC}$ contains a
unique graph, up to isomorphism, and this graph is then a ``canonical"
universal $\CC$-free graph.  Using well established model theoretic
terminology, we refer to this as the $\aleph_0$-categorical case.  The
role of ${\EE}_{\CC}$ in connection with the problem of determining
whether a universal $\CC$-free graph exists is explored in $\S 2$,
which
introduces the important technical notion of an {\sl existential type}.  
In $\S 3$ we characterize the $\aleph_0$-categorical case in terms of
the behavior of the associated algebraic closure operator on 
${\EE}_{\CC}$.  We begin the analysis of the algebraic closure
operator in
$\S 4$.  More delicate techniques for analyzing this operator are left to
the second part, as needed for applications.

Our thanks go to P.~Komjath for a close reading of a draft of the present
paper. 

\section 1. Existentially complete $\CC$-free graphs.

First we introduce some definitions and
notations which will be used in the whole paper.

\proclaim{Definition 1}
Let $\CC$ be a set of finite graphs.

1. A graph $G$ {\sl omits} $\CC$ if no subgraph of $G$ is isomorphic to
any graph in $\CC$.

2. $\GG_{\CC}$ is the class of all countable graphs omitting 
$\CC$.

3. A graph $G\in \GG_{\CC}$ is {\sl universal } (for
$\GG_{\CC}$) if every graph in $\GG_{\CC}$ is isomorphic
to an induced subgraph of $G$.
\endproclaim

\noindent{\bf Remarks}

1. There are two notions of universality which are generally considered.
We say that $G\in \GG_{\CC}$ is {\sl weakly universal} if every
graph in $\GG_{\CC}$ is isomorphic to a subgraph of $G$.  In
practice the two notions of universality behave similarly.  A universal
graph is evidently
weakly universal, and in practice proofs of the nonexistence of a
universal graph can often be doctored in standard ways to exclude weakly
universal graphs as well. 

For a theoretical analysis our definition of universality is to be
preferred, at least initially, as it facilitates the application of
general
methods.  To pass to the weakly universal case on a theoretical level is
in part a matter of replacing ``existential type" in $\S 2$ by ``positive
existential types", but the more pragmatic alternative of working mainly
with
(strictly) universal graphs on a theoretical level and then doctoring
specific construction is probably to be preferred.  

2. Let $T_{\CC}$ be the {\sl first order theory} of $\GG_{\cal
C}$.
Then the models of $T_{\CC}$ are all the $\CC$-free graphs and 
$\GG_{\CC}$ consists of the countable models of $T_{\CC}$,
which is a {\sl universal theory}.  This reflects the assumption that
all
graphs in $\CC$ are finite, and allows the application of model
theoretic methods. 

\proclaim{\bf Definition 2}
Let $\CC$ be a set of finite graphs.

1. If $G\subseteq H$ are graphs, we say that $G$ is {\sl existentially
complete} in $H$ if every existential statement $\phi$ which is defined in
$G$ and true in $H$ is also true in $G$.  Equivalently, if $A\subseteq B$
are finite induced subgraphs of $G$ and $H$ respectively, then there is an
embedding $f:B\rightarrow G$ taking $B$ isomorphically onto an induced
subgraph of $G$, with $f\restrictedto A$ the identity.

2. $G\in \GG_{\CC}$ is said to be {\sl existentially complete} 
(for $\GG_{\CC}$) if $G$ is existentially complete in each graph
$H$ such that $G\subseteq H\in \GG_{\CC}$. 

3. ${\EE}_{\CC}$ is the class of all existentially complete graphs
in $\GG_{\CC}$.

4. $T^{*}_{\CC}$ is the theory of ${\EE}_{\CC}$.  (In the proof
of Theorem 1 below we will determine this theory fairly precisely.)
\endproclaim

\proclaim{Example 1}
If $\CC=\emptyset$, then $\GG_{\CC}$ is the class of all
countable graphs and ${\EE}_{\CC}$ contains only one element up to
isomorphism: the random countable graph $G_{\infty}$ [Ra].  
$T_{\CC}$
is the theory of graphs, and 
$T^{*}_{\CC}$ is the theory of
$G_{\infty}$ (a complete theory).
\endproclaim

\proclaim{Example 2}
If ${\CC}=\{K_{3}\}$, a complete graph, then $\GG_{\CC}$ is the
class of countable triangle-free graphs and ${\EE}_{\CC}$ contains a
unique element up to isomorphism, called the generic triangle-free graph
$G_3$. $T_{\CC}$ is the theory of triangle-free graphs and $T^{*}_{\cal
C}$ is the theory of $G_3$.  For ${\CC}=\{K_{n}\}$, any $n$, the
situation is similar.
\endproclaim

\proclaim{Example 3}
If ${\CC}=\{K_{2}+K_{2}\}$, the disjoint sum of two copies of $K_2$,
then 
${\EE}_{\CC}$ contains two elements up to isomorphism: the triangle
$K_3$ and the star $S_{\infty}$ of infinite degree. The theory
$T^{*}_{\CC}$ is not a complete theory, since $K_{3}$ and 
$S_\infty$ have different theories.
\endproclaim

\proclaim{Example 4}
If ${\CC}=\{S_{3}\}$ ($S_n$ denotes a star of degree $n$ or order
$n+1$), then $T_{\CC}$ is the theory of graphs $G$ with vertex degree
at most 2, and 
$T^{*}_{\CC}$ is the theory of graphs in which every vertex has
degree 2, and which contain infinitely many cycles $C_n$ for each 
$n\ge 3$. The countable models $G$ of
$T^{*}_{\CC}$ are characterized up to isomorphism by the number of
connected components in $G$ isomorphic to a 2-way infinite path.  If $G_k$
is the model of $T^{*}_{\CC}$ with $k$ components of this form ($k\geq
0$),
then $G_{\infty}$ is universal for this class.   
\endproclaim

\proclaim{Remarks}
\rm

1. We will see below that $T^{*}_{\CC}$ is complete if the graphs in
$\CC$ are connected. This is the case of primary interest here.

2. It is easy to see that there is a universal graph in 
$\GG_{\CC}$ if and only if there is a universal graph in 
${\EE}_{\CC}$.  We will attempt to make this observation more useful
by analyzing $T^{*}_{\CC}$ and ${\EE}_{\CC}$ more clearly below.

3. The notion of existential completeness makes sense in almost any
context (though our reformulation in terms of embeddings is not always
accurate).  For example, algebraically closed fields are existentially
complete by Hilbert's Nullstellensatz; real closed fields are
existentially complete in the category of ordered fields (Tarski); and
dense linear orders are existentially complete in the category of linear
orders.

4. While it is natural to think of existential completeness as a form of
``algebraic closure", it does not involve the sort of {\sl finiteness
assumptions} connected intuitively with notion of algebraicity.  We will
introduce the model theoretic algebraic closure operator below.
\endproclaim

\proclaim{Theorem 1}
Let $\CC$ be a finite set of finite graphs.  Then 
\item{1.} ${\EE}_{\CC}$ is the class of countable models of the theory 
$T^{*}_{\CC}$.
\item{2.} If every $C\in {\CC}$ is connected, then $T^{*}_{\CC}$ is
complete.
\endproclaim

The proof will involve the general theory of model complete theories and
existentially complete structures, as presented in [HW].  We first give an
example showing the necessity of taking $\CC$ finite.

\proclaim{Example 5}
Let ${\CC}=\{C_{n}: n\ge 3\}$, the class of all cycles. Then ${\cal
G}_{\CC}$ is the class of countable forests and ${\EE}_{\CC}$ 
contains a unique graph $T_{\infty}$, up to isomorphism, a tree in which
every vertex is of countable infinite degree.  The models of $T^{*}_{\cal
C}$ are disjoint unions of any number of copies of $T_{\infty}$.
\endproclaim

\proclaim{Remark}
In Theorem 1, clause (2) follows from clause (1).  This is because clause
(1) provides one of the standard criteria for the theory $T^{*}_{\CC}$
to be model complete (Robinson's Test, [HW, Theorem 2.2]) and for such
theories, completeness is equivalent to the {\sl joint embedding
property}: any
two models of a theory should be contained as induced subgraphs in a third
[HW, Proposition 2.8].
If $\CC$ consists of connected graphs, then $T_{\CC}$ is closed
under the formation of disjoint sums. However connectedness is not a 
necessary condition for joint embedding:
\endproclaim

\proclaim{Example 6}
If ${\CC}=\{K_{3},K_{2}+K_{2}\}$ then $T_{\CC}$ has the joint
embedding property.  
\endproclaim

It is not clear whether one can easily recognize the
finite sets $\CC$ for which $T_{\CC}$ has the joint embedding
property.

The proof of Theorem 1 requires the following technical lemma.
Recall that a quantifier-free formula is {\sl conjunctive} if it is a
conjunction of atomic formulas and the negations of atomic formulas.
An existential formula of the form $\exists \bar{x} \phi$ with $\phi$ 
quantifier-free and conjunctive is called {\sl primitive}.  A typical
example of a conjunctive formula is a description of the isomorphism type
of an induced subgraph.

\noindent{\bf Notation.}

Let $\phi$ be a formula.  We write $T_{\CC}\vdash \phi$ if every 
$\CC$-free graph satisfies 
$``\forall \bar{x}\phi (\bar{x})"$ (we quantify over all free variables in
$\phi$).  In other words, $\phi$ is ``always" true in $\CC$-free
graphs.

\proclaim{Lemma 1}
Let $\CC$ be a finite set of forbidden substructures.  For each 
$n\ge 0$ there is a natural number $b_{n}$ such that for any two primitive
existential formulas $\phi,\psi$ such that 
\item{i.} $\phi$ contains at most $n$ existential quantifiers,
\item{ii.} $T_{\CC} \vdash\neg (\phi\wedge \psi)$, and  
\item{iii.} for each pair of variables $y_{1},y_{2}$
occurring in $\psi$, with at least one of them quantified, the clause
$y_{1}\not=y_{2}$ occurs as a conjunct in $\psi$,

\noindent 
there is a subformula $\psi_{1}$ of $\psi$ such that 
\item{1.} $\psi_{1}$ contains at most $b_{n}$ existential quantifiers.
\item{2.} $T_{\CC}\vdash\neg (\phi \wedge \psi_{1})$.
\endproclaim

We will first explain how Theorem 1 follows from this lemma, then prove
the lemma.  The following is essentially a corollary to Lemma 1.

\proclaim{Lemma 2}
Let $\phi(\bar x)$ be a universal formula.  Then there is an existential
formula $\psi(\bar x)$ such that
$$T_{\CC}\vdash\forall \bar{x}
[\phi(\bar{x})\longleftrightarrow\psi(\bar{x}) ].$$  
\endproclaim

\proof:
We use Proposition 1.6 (iii) of [HW].  Let $\Phi$ be the set of all
existential formulas $\phi '(\bar{x})$ such that
$$T_{\CC}\vdash \forall \bar{x}[\phi '(\bar{x})\longrightarrow
\phi(\bar{x})].$$
Then for 
$G\in{\EE}_{\CC}$, $\bar{u}\in G$, we have 
$$G\models \phi(\bar{u})\Longleftrightarrow G\models \phi '(\bar{u})
\hbox{ for some }\phi '\in \Phi.$$ 
In other words,
$$G\models \forall \bar{x}[\phi(\bar{x})\longleftrightarrow
\bigvee_{\Phi}\phi '(\bar{x})]. \eqno (*)$$
Note that the disjunction on the right is infinite; using Lemma 1 we will
replace $\Phi$ by a finite subset $\Phi '$ for which the analog of (*)
holds.
Thus with 
$\psi=\bigvee_{\Phi '}\phi '$, the claim follows.

Any existential formula is equivalent to a disjunction of primitive
existential formulas; so we may take $\Phi$ to consist of primitive
existential formulas. Similarly, the universal formula $\phi$ is
equivalent to 
a conjunction of negations of primitive existential formulas, so it
suffices to deal with the case $\phi=\neg \phi_{1}$ with $\phi_{1}$
primitive existential.  Finally, we may suppose that for each
$\phi'\in\Phi$ and each pair $y_{1},y_{2}$ of variables occurring
existentially quantified in $\phi'$, we have $y_{i}\not=y_{j}$
as a conjunct in $\phi'$ for $i\neq j$. Indeed if $\phi'=\exists
\bar{y}
\phi'_0(\bar{x},\bar{y})$, then 
$$\phi'\longleftrightarrow\exists \bar{y}(\phi'_0\wedge y_{i}=y_{j})
\vee \exists \bar{y}(\phi'_0\wedge y_{i}\not=y_{j})$$
so we may replace $\phi'$ if necessary by two disjuncts on the right and
then contract variables in the first disjunct.

After these preparations, $\Phi$ consists of formulas $\phi'$ to which
Lemma 1 applies, with $n$ the number of quantifiers occurring in $\phi$.
Thus if $\Phi'\subseteq \Phi$ consists of the primitive existential
formulas $\phi'$, 
in at most $b_{n}$ variables such that 
$$T_{\CC}\models\neg(\phi_{1}\wedge \phi'),$$
then
$$T_{\CC}\models
\forall\bar{x}[\phi(\bar{x})\longleftrightarrow
\bigvee_{\Phi'}\phi'(\bar{x})].$$ 
\qed

\proof\ of Theorem 1:
By Lemma 2, every universal formula is equivalent to an existential
formula modulo $T_{\CC}$. This is equivalent to clause (1) of Theorem
1 by [CK, Theorem 3.5.1].  As noted before, clause (2) follows from clause
\quad (1).
\qed

\proof\  of Lemma 1:
We proceed by induction on $n$, the number of quantified variables in
$\phi(\bar{x})$.  Let $k=max\{|C|:C\in{\CC}\}.$

If $n=0$, then $\phi$ is quantifier free and we will take $b_0=k$.
Suppose $T_{\CC}\models \neg (\phi\wedge\psi)$.
We have $\psi=\exists \bar{y}\psi_0(\bar{x},\bar{y})$,
$\psi_0$ either contradicts $\phi$ explicitly, or states
that the induced graph on some subset of $k$ vertices from
$\bar{x},\bar{y}$ contains a forbidden subgraph.  In the former case
$\psi$ can be replaced by a quantifier free formula, and in the latter
case by a formula in at most $k$ quantified variables.

For the induction step, let
$\phi=``\exists\bar{y}\phi_0(\bar{x},\bar{y})"$ have $n+1$ quantified
variables and let $\psi = ``\exists\psi_0(\bar{x},\bar{y'})"$.
Let $A$ and $B$ be the graphs on vertices $\bar{x},\bar{y'}$ described by
$\phi_0$ abd $\psi_0$ respectively, that is, edges exist as specified
by $\phi_0$ and $\psi_0$. 
As $T_{\CC}\models
\neg(\phi\wedge\psi)$ the free joint of $A$ and $B$ over $\bar{x}$
contains a forbidden graph $C\in {\CC}$. For each pair of variables
$y_{i}$ in $C\cap\bar{y}$ and $y_{j}'$ in $C\cap \bar{y}'$, introduce a
new variable $x_{ij}$ and let $\phi^{*}_0(\bar{x},x_{ij},\hat y)$ and 
$\psi^{*}_0(\bar{x},x_{ij},\hat y')$ be obtained by replacing $y_{i}$
by
$x_{ij}$ in $\phi_0$ and $y_{j}'$ by $x_{ij}$ in $\psi_0$.  Thus
$\hat y$
and $\hat y'$ are $\bar{y}$ and $\bar{y}'$ with $y_{i}$ or $y_{j}'$
deleted.
Write $\hat x$ for $\bar{x},x_{ij}$.

Let $\phi^{*}=\exists \hat y \phi^{*}_0(\hat x,\hat y)$ and 
$\psi^{*}=\exists \hat y'\psi^{*}_0(\hat x,\hat y')$.
Then $\phi^{*}$ has $n$ quantified variables and 
$T_{\CC}\models \neg(\phi^{*}\wedge\psi^{*})$, since any model of
$T_{\CC}\cup\{\phi^{*},\psi^{*}\}$ gives rise to a model of 
$T_{\CC}\cup\{\phi,\psi\}$; the variables $y_{i},y_{j}'$ may be
realized by the value of $x_{ij}$.

By induction hypothesis for each choice of $i$ and $j$, $\psi^{*}$
contains a subformula $\psi^{*}_{ij}$ involving at most 
$b_n$ variables
so that $T_{\CC}\models \neg(\phi^{*}\wedge\psi^{*}_{ij})$.

Let $\bar{y}''\subseteq \bar{y}'$ be the set of at most $k+k^{2}b_{n}$
variables
consisting of $C\cap \bar{y}'$ together with the all quantified
variables
from any $\psi^{*}_{ij}$, and let $\psi_{1}$ be the restriction of $\psi$
to $\bar{y}''$.  Then we claim
$$T_{\CC}\models \neg(\phi\wedge\psi_{1}) \eqno (*)$$   
so we may take $b_{n+1}=k+k^{2}b_{n}$.

For (*), consider any model $\cal M$ of $\phi\wedge\psi_{1}$. Then $C$
embeds in the free join over $\bar{x}$ of the
induced graphs $A,B$ on $\bar{x},\bar{y}$ and $\bar{x},\bar{y}''$. 
So if $\cal
M$ omits $\CC$, there must be some identification $x_{i}=y_{j}$ with
$x_{i},y_{j}\in C$. This is exactly what is ruled out by $\psi^{*}_{ij}$.
\qed

\proclaim{Corollary to Theorem 1}
Let $\CC$ be a finite class of finite graphs.  Then $T^{*}_{\CC}$ is
model complete and is the model companion of $T_{\CC}$.
\endproclaim

\proof:
This is equivalent to Theorem 1, part (1).
\qed

\section 2. Universal Graphs and existential types.

In this section we give criteria for the existence of a universal graph in 
$\GG_{\CC}$, for $\CC$ a finite set of finite connected
graphs.  We will show that when there is a universal graph in 
$\GG_{\CC}$, there is a {\sl canonical} one, namely the
``$\aleph_0$-saturated" graph in ${\EE}_{\CC}$.  We will also show the
relationship of this problem to a model theoretic notion of {\sl algebraic
closure}. We review the definitions.

\proclaim{Definition 3}
\rm

Let $\CC$ be a collection of finite forbidden subgraphs.

1. The {\sl existential type} $tp_{G}(\bar{a})$ of a finite sequence
$\bar{a}=a_{1},a_{2},\cdots, a_{n}$ in a graph $G\in {\EE}_{\CC}$ is
the set of existential formulas $\phi(\bar{x})$ such that
$G\models\phi(\bar{a}).$ The Stone space $S_{n}(T^{*}_{\CC})$ is the
set of all existential types $tp(\bar{a})$ of sequences
$\bar{a}=a_{1},\cdots,a_{n}$ 
in any graph $G\in {\EE}_{\CC}$.

2. $G\in{\EE}_{\CC}$ is $\aleph_0$-{\sl saturated} if for all $n$, all
$\bar{a}\in G$ of length $n$, and all $(n+1)$-types in
$S_{n+1}(T^{*}_{\CC})$ whose restriction to the first $n$ variables is
$tp_{G}(\bar{a})$, there is $v\in V(G)$ so that $tp_{G}(\bar{a},v)$ is the
specified type.
\endproclaim

\proclaim{\bf Example 7}
When ${\CC}=\{S_{3}\}$, specifying the type of an element $a$
in $G\in {\EE}_{\CC}$ is equivalent to 
describing the isomorphism type of its
connected component in $G$. In particular if $a_{1},\cdots,a_{n}$ lie in
distinct components isomorphic to 2-way infinite paths,
$\omega$-saturation yields an element $a_{n+1}$ lying in another
component isomorphic to such a path.  Thus the $\aleph_0$-saturated model is
the largest model in ${\EE}_{\CC}$.  This is the case in general. 
\endproclaim

\proclaim{Theorem 2}
Let $\CC$ be a finite set of connected forbidden subgraphs.  Then the
following are equivalent:
\endproclaim

1). There is a universal graph in $\GG_{\CC}$.

2). There is a universal graph in ${\EE}_{\CC}$.

3). ${\EE}_{\CC}$ contains a unique $\aleph_0$-saturated graph, up to
isomorphism.

4). $S_{n}(T^{*}_{\CC})$ is countable, for any $n$.

\proof:
This is a special case of general model theoretic principles [CK,
$\S 2.3$]. We sketch the ideas here.

The equivalence of 1) and 2) is immediate. It suffices to note that any 
$G\in \GG_{\CC}$ embeds into a $G^{*}\in {\EE}_{\CC}$.
For the equivalence of 2) to 4) one recalls that ${\EE}_{\CC}$ is
the class of countable models of $T^{*}_{\CC}$. We will show 
$2)\Rightarrow 4)\Rightarrow 3)\Rightarrow 2).$

$2)\Rightarrow 4)$. Let $G\in {\EE}_{\CC}$ be universal.  As $G$ is
countable, the set $\{tp_{G}(\bar{a}):\bar{a}=a_{1},\cdots,a_{n}\in G\}$
is countable.  Any type $tp_{G'}(\bar{a})$ realized in any 
$G'\in {\EE}_{\CC}$ will be realized in $G$ since we may take $G'$
to be
an induced subgraph of $G$ by universality and
$tp_{G}(\bar{a})=tp_{G'}(\bar{a})$ by existential completeness.

$4)\Rightarrow 3)$.
If $S_{n}(T^{*}_{\CC})$ is countable for all $n$, one builds a
countable saturated model as the limit of an 
increasing countable sequence of models in
${\EE}_{\CC}$, see [CK, Theorem 2.3.7].

The uniqueness follows from the completeness of $T^{*}_{\CC}$ 
[CK, Theorem 2.3.7].

$3)\Rightarrow 2)$.
Saturated models are universal [CK, Theorem 2.3.10].
\qed

In the examples, one often encounters the special case in which 
${\EE}_{\CC}$ contains a unique model up to isomorphism, so that the
$\aleph_0$-saturation condition is vacuous.  This is a rather special case 
in model theory, and the frequency of its occurrence in our context is an 
indication that something more specialized is involved. To analyze this
further we introduce the notion of {\sl algebraic closure}.

\proclaim{Definition 4}\ 

Let $\CC$ be a set of forbidden subgraphs, $G\in{\EE}_{\CC}$,
$A\subseteq G$, $a\in G$.  We say that $a$ is {\sl algebraic} over $A$ (in
$G$)
if there is an existential formula $\phi(x,\bar{a})$ with $\bar{a}\in A$
such that the set $\{a'\in G:\phi(a',\bar{a})\}$ is finite and contains
$a$. We write 
$\acl_{G}(A)$ (algebraic closure) for the set of 
$a\in G$ algebraic over $A$.  We say $A$ is 
{\sl algebraically closed in G} if $\acl_{G}(A)=A$.
\endproclaim

\proclaim{Lemma 3}
Let $\CC$ be a finite set of connected forbidden subgraphs.  If ${\cal
G}_{\CC}$ contains a universal graph then the set of isomorphism types
of induced subgraphs of graph $G\in {\EE}_{\CC}$ on subsets of the
form $\acl(A)$ with $A$ finite, is countable.
\endproclaim

\proof:
Let $G\in {\EE}_{\CC}$ be universal. Then for any $G'\in {\cal
E}_{\CC}$ and any $A\subseteq G'$ finite, an embedding $\iota$ of $G'$
into $G$ given an isomorphism between 
$G'\restrictedto \acl_{G'}(A)$ and 
$G\restrictedto \acl_{G}(\iota A)$.  The
point here is that $\acl_{G}(\iota A)=\iota[\acl_{G'}(A)]$, by existential
completeness.
\qed
    
It would be pleasant if the converse held: in other words, to show the
nonexistence of universal graphs one would be obligated to construct
uncountably many isomorphism types of algebraic closures of finite sets.
This is what has actually occurred in all examples treated to date
[ChK,CS2,CST,FK1,Ko,FK2,GK,KP1].

In fact, in most cases one of the two extremes of the following
pseudo-dichotomy have been encourtered:

I. The algebraic closure of a finite set is finite.

II. There are uncountably many isomorphism types of induced subgraphs on
sets $\acl_{G}(A)$, with $A$ finite, in graphs $G\in {\EE}_{\CC}$.

On the other hand the example ${\CC}=\{S_{3}\}$, a star of degree 3,
shows that case I is indeed a special case, as one might anticipate. This
makes it all
the more surprising that this case is typical in practice, in
contexts where universal graphs exist.

All of this leaves open the possibility, already referred to, that 
case II is an
exact criterion for the nonexistence of universal graphs.  To refute this
in
the category of graphs is not so easy.  We will give an example in the
category of vertex-colored graphs.  It should not be too difficult to
encode this
as an example in the category of graphs, but it would be more to the point
to prove the general encoding conjecture noted in the introduction, which
we will not undertake here.

\proclaim{Example 8}
We work with vertex colored graphs in which there are three colors: 0, 
+1, -1. Each vertex of color 0 has at most two neighbors of color 0, and
only one of
the other two colors occurs among its neighbors.  Vertices of colors +1
and
-1 are adjacent to at most one vertex, which must have color 0.  This
clearly corresponds to a finite set of connected forbidden subgraphs.  In
${\EE}_{\CC}$ the graphs consist of cycles and 2-way infinite paths
made up of vertices of color 0, each decorated with infinitely many
adjacent
vertices of color +1 or -1.  It is easy to see that the algebraic closure
of a finite set $A$ consists of the union of the connected components
of vertices of color 0 in or adjacent to $A$, together with vertices in
$A$
of color +1 and -1.  Thus there are countably many induced subgraphs on
$\acl(A)$ for $A$ finite.

However, the type of an element $v$ of color 0 contains a specification
of the colors +1 and -1 of the neighbors of all vertices of color 0 
in its connected component.  Thus $S_{1}(T^{*}_{\CC})$ is uncountable.
It follows that the types in general contain information not controlled by
the algebraic closure operation.
\endproclaim

On the other hand, we will show that when condition (I) holds, control of
algebraic closure {\sl is} enough.
Indeed, in the example just discussed, there are only countably many types
associated with vertices of color 0 whose connected component, among the
vertices of type 0, is finite.  In fact, if the order of the connected
component in question is specified, there are finitely many possible
types.

\section
3. $\aleph_0$-categoricity and local finiteness

A theory is said to be $\aleph_0$-{\sl categorical} if it has a unique
countable model, up to isomorphism.  As we have noted, among theories of
the form $T^{*}_{\CC}$ for which a universal countable model exists,
the $\aleph_0$-categorical case is surprisingly common.  The next result 
casts some light on this phenomenon.

\proclaim{Theorem 3}
Let $\CC$ be a finite set of connected finite graphs.  Then the
following are equivalent:
\endproclaim

(1). $T^{*}_{\CC}$ is $\aleph_0$-categorical.

(2). $S_{n}(T^{*}_{\CC})$ is finite for each $n$.

(3). For $A\subseteq {\cal M}\models T^{*}_{\CC}$ finite, we have
$\acl(A)$ finite.

\noindent{These conditions imply}

(4). $\GG_{\CC}$ contains a universal countable graph.

By Theorem 1, (2), $T^{*}_{\CC}$ is complete.  Therefore the
equivalence of (1) and (2) holds by general model theory [CK,Theorem
2.3.13.]. That
(1) implies (4), and (2) implies (3), are both immediate.  Thus all that
requires proof is the implication from (3) to (2).  For this we  prove a
more refined technical lemma, based on  the following definition and fact.

\def\tp{\hbox{\rm tp}}
\proclaim{Definition 5}
Let $G$ be a graph, and $A\subseteq G$. Set
$$\eqalign{
\tp_{n}(A)=\{\phi(\bar{a}):
	&\hbox{ $\phi$ is existential, with at most $n$ quantified
	 variables,}\cr 
	&\hbox{$\bar{a}\in A$, and $\phi(\bar{a})$ holds in $G$} \}.\cr
}$$
(This depends on $G$, and one may write $\tp^{G}_{n}(A)$ to show this
dependence.)  
\endproclaim

\proclaim{Fact 1. (Park, cited in [Ba])}
Let $A$ be algebraically closed in $B$.  Then there is 
$C\succ B$ and $B'\simeq B$ (over $A$) with $B'\prec C$ and $A=B\cap B'$.
Note that in [Ba] the term ``Park-a.c." is used for our ``algebraically 
closed".
\endproclaim

\proclaim{Lemma 4}
Let $\CC$ be a finite set of finite graphs, and 
$A\subseteq G\in {\EE}_{\CC}$ 
with $A$ finite and algebraically closed.  Then for
$n=\max\{|C|:C\in {\CC}\}$, $\tp_{n}^{G}(A)$ determines $\tp(A)$.
\endproclaim

\proof:
We write $\bar{a}$ for $A$ arranged as a finite sequence. Let
$e(\bar{a})$ be an existential sentence.
We claim that $e(\bar{a})$ holds in $G$ if
and only if the following theory $T_{e}$ is consistent:
$$``A \hbox{ is algebraically
closed"}\cup T_{\CC}\cup \tp^{G}_{n}(\bar{a})\cup
\{e(\bar{a})\}.\leqno (T_{e})$$
One may easily find axioms expressing
the assertion that $A$ is algebraically closed. Thus $T_{e}$ is
indeed a
first order theory.  If $e(\bar{a})$ holds in $G$, then
$T_{e}$ holds in $G$ and thus $T_{\CC}$ is consistent.

Suppose conversely that
$T_{e}$ holds in some $G_{1}$. We claim that $e(\bar{a})$ will then hold
in $G$.

Let $e(\bar{a})=\exists\bar{y}e_0(\bar{a},\bar{y})$ with $e_0$
quantifier-free.  We may suppose that $e$ is primitive, and $e_0$ is
conjunctive.  
Choose $\bar{b}$ in $G_{1}$ so that 
$e_0(\bar{a},\bar{b})$ holds.  We may suppose 
$\bar{b}\cap A=\emptyset$, adjusting $e_0$ if necessary. 
Form $G'=G\cup \bar{b'}$ by freely amalgamating $G$ with a copy 
$\bar{a}\bar{b'}$ of
$\bar{a}\bar{b}$ over $\bar{a}$. That is, the
edges in $G'$ lie in $G$ and in $\bar{a}\bar{b'}$.  Note that $G$ and
$G_{1}$ agree on $\bar{a}$, as a description of the induced graph on
$\bar{a}$ is contained in $\tp_{n}(\bar{a})$. 

If $G'\in \GG_{\CC}$ then as $G\subseteq G'$, 
$G\in {\EE}_{\CC}$, and $e(\bar{a})$ holds in $G'$ we find that 
$e(\bar{a})$ holds in $G$, as claimed. Suppose now that 
$G'\not\in \GG_{\CC}$.  Then we will show that 
$G_{1}\not\in \GG_{\CC}$, contradicting our assumption on $G_{1}$.

We have some $C\in {\CC}$ which embeds into $G'$, and we may take
$C\subseteq G'$.  Let 
$\bar{c_0'}=C-(A\cup \bar{b}')\subseteq G$, and let 
$\phi_0(\bar{a},\bar{c_0}')$ be a conjunctive quantifier-free formula
specifying the isomorphism type of the induced subgraph on
$\bar{a},\bar{c_0}'$.  Then the existential formula 
$$\phi(\bar{x})=``\exists \bar{y}\phi_0(\bar{a},\bar{y})"$$ 
belongs to 
$\tp_{n}^{G}(A)$. Hence we have 
$\bar{c_0}$ in $G_{1}$ satisfying 
$\phi_0(\bar{a},\bar{c_0})$.

As $\bar{c_0}\cap A=\emptyset$ and $A$ is algebraically closed in
$G_{1}$,
by repeated applications of Fact 1 we can find disjoint sequences 
$\bar{c}_0^{(1)},\cdots, \bar{c}_0^{(k)}$ in $G_{1}$, for any $k$, so
that the induced subgraphs on $\bar{a}\bar{c}_0^{(i)}$ are isomorphic to
$\bar{a}\bar{c_0}$ in the natural order.

Choose $k>|\bar{b}|.$ Then for some $i$, $\bar{c}_0^{(i)}\cap
\bar{b}=\emptyset $, and thus the free amalgam of $\bar{a}\bar{c}_0$
with
$\bar{a}\bar{b}$ over $\bar{a}$ embeds into the induced graph on
$\bar{a}\bar{b}\bar{c}_0^{(i)}$.  But this free amalgam is also 
isomorphic to the subgraph of $G'$ induced
on $\bar{a}\bar{b}'\bar{c_0}'$, which is $C$.  Thus 
$C$ embeds in $G_{1}$, a contradiction.
\qed

\proof\  of Theorem 3:
As noted above, we need only check $(3)\Rightarrow (2)$.  Assuming
(3), then for $n$ fixed there is a uniform bound on $|A|$ for $A$ the
algebraic closure of a set of $n$ elements in a model of $T^{*}_{\CC}$.
Thus it suffices to show that for each such $A$, the type of $A$ in a graph
$G\in {\EE}_{\CC}$ is determined up to  finitely many possibilities.
Indeed, with $A$ fixed, by the preceding lemma there is  $N$ such that 
$\tp_{N}(A)$ determines $\tp(A)$; and there are only finitely many
possibilities for $\tp_{N}(A)$.
\qed

Thus if the algebraic closure operation is uniformly locally finite on
${\EE}_{\CC}$, a universal graph exists.  Earlier we showed by
example that when it is not uniformly locally finite, knowledge of this
operator does not in general settle the question of the existence of a
universal graph: of course, at the other extreme (case II of \S 2),
the question is also settled by the structure of algebraic closure.

\section 4. Algebraic closure.

In view of the importance of the algebraic closure operator in dealing
with problems of universality, it is worth while making explicit what is
involved.

\proclaim{Definition 6}
Let $A,B$ be graphs and $f:V(A)\longrightarrow V(B)$.
Then $f$ is a {\sl homomorphism} if $f$ carries edges to edges.
\endproclaim

\noindent{\bf Remarks}

1.  An injective homomorphism is an isomorphism with a subgraph (not
necessarily induced).

2. We deal throughout with graphs without loops.  In particular if a
homomorphism $f:A\longrightarrow B$ identifies two vertices of $A$, they
cannot be linked by an edge.  (We could just as well allow loops. In this
case, if the loop on one vertex is in $\CC$, we recover the loop-free
context.)

\proclaim{Lemma 5}
Let $\CC$ be a finite collection of finite graphs, and 
$A\subseteq G\in {\EE}_{\CC}$. Then the following are equivalent:
\endproclaim

(1). $A$ is not algebraically closed in $G$.

(2). There is some $C\in {\CC}$ and a homomorphism
$C\longrightarrow C'\subseteq G$ so that $C$ embeds in the free
amalgam over
$A$ of $|C|$ copies of $C'$.

\proof:
$(2)\Longrightarrow (1):$ Let 
$h:C\longrightarrow C'$ as in (2) and let $B=C'-A$. If $G-A$ contains
$|C|$ disjoint copies $B^{i}$ of $B$ (isomorphic over $A$), then the free 
amalgam of
$|C|$ copies of $C'$ over $A$ embeds into $A\cup\bigcup_{i\leq|C|} B^{i}$,
and hence $C$ embeds in $G$, a contradiction.  By Park's Theorem (Fact 1),
our claim follows.

$(1)\Longrightarrow (2)$:
As $A$ is not algebraically closed, there is $b\in \acl(A)-A$, and there is
an existential formula 
$\phi(\bar{a},b)=``\exists \bar{y}\phi_0(\bar{a},b,\bar{y})"$ so that 
$|\{b'\in G: \phi(\bar{a},b')\}|=k<\infty.$ Let 
$\bar{b}\in G$ satisfying $\phi_0(\bar{a},b,\bar{b})$, and set 
$B=\{b\}\cup \bar{b}$. With a slight change of notation, we may suppose 
$B\cap A=\emptyset$.

Let $G_0=AB^{1}\cdots B^{k+1}$ be the free amalgam over 
$A$ of $k+1$ copies
$AB^{i}$ of $AB$ (isomorphic over $A$).  Let $G_{1}$ be the free amalgam 
of $G$ and
$G_0$ over $A$.  Then $G_{1}\not\in \GG_{\CC}$, as otherwise
after
extending to $G_{2}\in \GG_{\CC}$, we find $G\prec G_{2}$ but 
$|\{b'\in G_{2}:\phi(\bar{a},b')\}|>k$, a contradiction.

As $G_{1}\not\in \GG_{\CC}$, there is $C\in{\CC}$ and an
embedding $f:C\longrightarrow G_{1}$. We alter this to a homomorphism
$h:C\longrightarrow G$ by mapping each $B^{i}$ isomorphically over $A$ to
$B$. Let $C'$ be the image of $h$.  Then the free join of $|C|$ copies of
$C'-A$ over $A$ contains the image of $f$, as required.
\qed

We give a simple example to illustrate the power of this result.  Later as
we go into applications in more detail, we will get considerably more 
mileage out of the same idea.

\proclaim{Theorem 4}
Let $\CC$ be a finite set of connected finite graphs.  Suppose that for
any $C\in{\CC}$ and any surjective homomorphism 
$h:C\longrightarrow C'$, that $C'$ contains a graph in $\CC$.  Then for 
$A\subseteq G\in {\EE}_{\CC}$, $\acl(A)=A$.  In particular,
$T^{*}_{\CC}$ is $\aleph_0$-categorical and hence there is a
universal graph in $\GG_{\CC}$.
\endproclaim

\proof:
If $A\subseteq G$ is not algebraically closed, application of Lemma 5
produces $h:C\longrightarrow C'\subseteq G$, but then 
$G\not\in \GG_{\CC}$.
\qed

\proclaim{Example 9}
Fix $k$.  Let $\CC$ consist of all cycles of odd lengths, up to $2k+1$.
Then there is a universal graph in $\GG_{\CC}$.
\endproclaim

This result was first proved in [KMP] with an elaborate
amalgamation argument, containing some minor
inaccuracies which were subsequently corrected.  This should serve to
illustrate the utility
of our general considerations.  We will use the same idea below to
construct a number of new examples.

\part II.  Applications

In this part, we first review the known results from our point of view.
>From this point of view, the main question is the behavior of the
algebraic closure operation on finite sets.  This qualitative problem can
be rephrased in quantitative terms; from that point of view, the known
results leave open a number of questions regarding the estimates for the
size of $\acl(A)$ in terms of $|A|$, and similar issues, which we will
point out in detail.

\section
5.  Negative results: explosion of algebraic closure.

The negative results all depend on the construction of $2^{\aleph_0}$
nonisomorphic induced graphs of the form $\acl(A)$ for $A$ of some fixed
size, which can be read off explicitly from the various papers, though the
terminology varies somewhat.  In such cases there are two natural
questions concerning $|\acl(A)|$:

(I).  What is the least cardinality $\alpha$ such that there are 
$2^{\aleph_0}$ possible isomorphism types for the graph induced on
$\acl(A)$ in a graph $G\in \EE_{\CC}$, with $|A|=\alpha $?

(II).  What is the least cardinality $\alpha'$ such that $\acl(A)$ is
infinite in some $G\in \EE_{\CC}$, with $|A|=\alpha'$?

One suspects these are usually equal, though exceptions were mentioned
earlier.  All of the negative results on universal graphs to date may be
phrased as explicit upper bounds on $\alpha$ in various cases.

\noindent{\bf 5.1.\sl  2-connected graphs.}

The main result of [FK1] gives a bound for $\alpha$ when 
${\CC} = \{C\}$ consists of a single constraint $C$ which is
2-connected
and not complete: 
$$ \alpha \le 4(2k)^{2g}. \leqno (1)$$ 
with $k=2|V(G)|-1, g=|V(G)|+1$.
Actually the result proved is significantly more general. The same bound
is given when $C$ contains a block $C_0$ which is 2-connected, and which
contains two nonadjacent vertices $u,v$ so that $C_0$ does not embed in 
$C_{uv}$, the graph obtained from $C$ by identifying $u$ with $v$ (keeping
all edges).

Various special cases proved earlier give sharper estimates for more
specific constraints.  In [KP1] one finds $\alpha \le 8m-7$ when
$C=K_{m,n}$ is complete bipartite $(m\le n)$.  In [ChK] one finds 
$\alpha \le 4N+1$ with $N=(14^{\nu}-1)/13, \nu = |V(G)|$, for $G$ a cycle
of
length at least 4; and the same bound is obtained in [CS2] when $\CC$
is a finite set of cycles, taking $\nu =\max\{|V(G)|:C\in {\CC}\}.$
There is one exception in this case: when $\CC$ consists of all odd
cycles up to some bound, there is a universal graph $(\alpha = \infty)$;
this was mentioned  above, following Theorem 4.

The special case in which $\CC$ consists of all cycles up to some even
bound was considered in [GK]; they found $\alpha\le 5$ in this case.

All of this raises a number of natural questions.  First of all, can one
combine [FK1] and [CS2] to identify all finite sets $\CC$ of
2-connected graphs for which there is a corresponding universal graph, and
to estimate $\alpha$ in the other cases?
   
Secondly, can one obtain a respectable lower bound for $\alpha$, or at
least determine whether $\alpha$ is unbounded in most cases?
Some information is provided by the following:

\proclaim{Lemma 6}
Let $\CC$ be a finite set of $k-$connected graphs.  Let 
$G\in {\EE_{\CC}}$ and $A\subseteq V(G), |A|<k$.  Then $\acl(A)=A$.
In particular $\alpha \ge k.$
\endproclaim

\proof:
We apply Lemma 5.  If $C\in {\CC}$ embeds in a free amalgam of copies
of $C'$ over $A$, then as $C$ is $k-$connected with $k>|A|, C$ would embed
in $C'$, hence in $G$.
\qed

\proclaim{\bf Example 10}
1.  If $\CC$ is a finite set of cycles, this tells us only that
$\acl(A)=A$ when $|A|=1$.

2.  If ${\CC} = \{K_{m,n}\}$ with $m\le n$ we find $\acl(A)=A$ for
$|A|<m$, and $\alpha \ge m$.  This matches the upper bound in [KP1] 
reasonably well.

3.  If ${\CC} =\{C\}$ with $C$ a complete graph $K_{n}$ with one edge
deleted, we find $\acl(A)=A$ for $|A|\le n-2$, so $\alpha \ge n-1$.
\endproclaim

This leaves a rather large gap between the upper and lower bounds for
$\alpha $ in most cases.  One suspects the upper bounds could be 
sharpened considerably.

\noindent{\bf 5.2. \sl Trees.}

A tree is called {\sl bushy} if it has no vertices of degree 2.  For
constraint sets $\CC$ consisting of a single bushy tree with at least
5 vertices, the result of
[CST] yields a bound slightly sharper than the following:
$$\alpha < n. \leqno (2)$$
This can be radically improved: if $\alpha<\infty$ then $\alpha=1$
in the case of trees (Proposition 6 below).

One peculiarity of tree constraints is that for any 
$G\in {\EE_{\CC}}$ (where $\CC$ 
consists of a single tree constraint $T$) we never have $\acl(A)=A$ when
$|A|=1$, unless $|V(T)|=2$.  This can be seen using Lemma 5.

\section
6.  Positive results: local finiteness estimates

In most of the cases in which a universal graph is known to exist,
$T^{*}_{\CC}$ is $\aleph_0$-categorical, and the situation is
described by Theorem 3.  In such cases the criterion in part (3) of that
theorem has not been used.  Indeed there are a variety of approaches to
$\aleph_0$-categoricity and it does not seem reasonable to insist on one
as
most appropriate in all cases, but we have indicated some situations in
which the computation of algebraic closure is effective, following Theorem
4.  We will give some new applications in the following section. 
Here we review the known positive results, with an eye on the additional
information they furnish about algebraic closure in such cases.  The
natural problem here is to estimate the function
$$c(n)=\max\{|\acl_{G}(A)|:|A|=n,A\subseteq V(G),G\in {\EE}_{\cal
C}\}.$$
Here upper bounds are the main point, but one may look for
accurate asymptotics.

\noindent{\bf 6.1.\sl Trees.}

Tallgren has conjectured that the only trees $T$ for which 
$\GG_{\{T\}}$ has a universal object are the paths and the trees
obtained from a path by attaching one additional edge.  His proof of the
existence of a universal graph in the latter case is unpublished, but this
case is of considerable interest as it affords an example in which 
$T^{*}_{\CC}$ is not $\aleph_0$-categorical, but a universal graph
exists.  This point is illustrated quite well by the simple example of a
star $S_{3}$ of degree 3, discussed as Example 4 in $\S 1$.

\noindent{\bf 6.2.\sl 2-connected case.}

Previously only two examples of finite sets $\CC$ of 2-connected graphs
were known for which $\GG_{\CC}$ has a universal graph: ${\cal
C}=\{K_n\}$, a single complete graph, or ${\CC}=\{C_{2k+1}:k\leq n\}$
the set of odd cycles of size up to to some bound.  Both are covered by
Theorem 4, as noted earlier: indeed $\acl(A)=A$ for all $A$, and $c(n)=n$.
Additional examples arising from Theorem 4 will be considered in the next
section.

\noindent{\bf 6.3.\sl  Bow-Ties}

Any graph can be analyzed as constructed from a tree of ``blocks"
(2-connected graphs and edges).  However, we know of no way to combine the
analysis of 2-connected constraints and tree constraints to produce
something more general.  For that matter, relatively few explicit examples
have been successfully analyzed to date.  Komj\'{a}th
[Ko]
did find one example in which a universal graph exists.  Such examples are
presumably quite rare.  We will give new examples in \S 8.  Here we give
an analysis of the F\"{u}redi-Komj\'{a}th ``bow-tie" example in terms of
our machinery of algebraic closure.

A {\sl bow-tie} is the graph on five vertices formed by attaching two
triangles
to a common vertex.  More generally, one may consider bouquets of complete
graphs with one common vertex.  For bouquets of at least three 
complete graphs of
constant size, it is shown in [Ko] that the only ones corresponding to
universal graphs are the bow-tie and the degenerate bouquets consisting of
one complete graph.

Let $B$ be the bow-tie, ${\CC}=\{B\}$.  We show that $T_{\CC}^{*}$ is 
$\aleph_0$-categorical, and in particular there is a universal countable
bow-tie-free graph.  This follows by combining Theorem 3 with the following
estimate.

\proclaim{Proposition 1}
Let $G\in {\EE}_{\{B\}},A\subseteq G$ finite.  Then 
$|\acl(A)|\leq 4|A|$.
\endproclaim

\proof:
Call an edge of $G$ {\sl special} if it lies in two
triangles of $G$.  We make the following claims, which will be verified
below.

(1). Every triangle in $G$ contains at least one special edge.  

(2). Every point that lies on a triangle, but no special edge of that
triangle, lies on a unique triangle.

(3). If a point lies on two special edges, it lies on a graph
$K\cong K_{4}$.  In this
case, any triangle containing that point is contained in $K$.

Assuming these claims for the moment, we proceed as follows.  Given    
$G\in{\EE}_{\{B\}}$, $A\subseteq V(G)$ finite, let $A^{*}$ be the union
of $A$ with the set of all vertices of $G$ which lie on special edges
which themselves lie on triangles containing a point of $A$.  It follows
from
(2),(3) that 
$$|A^{*}|\le 4|A|.\leqno(4)$$
Thus it will suffice to show that $A^{*}$ is algebraically closed.  We
show first $$A^{**}=A^{*}.\leqno(5)$$

Let $u\in A^{**}-A.$ Then $u$ lies on a special edge $e$, where $e$ lies on a
triangle $t$ meeting $A^{*}$.  Let the vertices of $e$ be $\{u,v\}$, and
let the third vertex of $t$ be $w$.

We claim $u\in A^{*}$.  Assume not.  Then $v$ or $w$ belongs to $A^{*}$,
but $t$ contains no vertex of $A$.  Thus $v$ or $w$ is in $A^{*}-A$.

If $v\in A^{*}-A$, then $v$ lies on a special edge $e'$ which lies on a
triangle $t'$ meeting $A$ in a vertex $a$.  If $e'=e$ this forces $u\in
A^{*}$, as desired.  If $e'\ne e$ then by (3) $v$ lies on a $K_{4}$,
containing $t$ and $t'$. 
Hence $a,u,v$ are the vertices of a triangle in
$G$ and therefore $u\in A^{*}$, as claimed.

If $w\in A^{*}-A$ and $w$ 
lies on a special edge of $t$, then by (2) all edges of $t$ are
special, and then the argument above applies to $w$.  If not, then by (2)
$w$ lies on a unique triangle.  Then if $w\in A^{*}$ then $t$ meets $A$, a
contradiction.

Now we show
$$A^{*}\hbox{ is algebraically closed. }\leqno(6)$$

We apply Lemma 5.  $B$ has only two proper homomorphic images, so
applying the criterion of Lemma 5, if $A^{*}$ is not algebraically
closed then there is a triangle $t$ meeting $A^{*}$ in one vertex.  By (1)
$t$ contains a special edge, and some vertex of that edge then lies in 
$A^{**}-A^{*}$, a contradiction.

It remains to verify our claims (1)-(3).  Both (2) and (3) are direct
consequences of the assumption that $G$ is $B$-free, by inspection.
We turn to (1).

Let $e$ be an edge of a triangle $t$ lying in $G$. Let $G^{*}$ be the
graph formed from $G$ by attaching an additional triangle containing the
edge $e$.  In $G^{*}$, $e$ is special.  If $G^{*}$ is $B$-free, then as
$G\in {\EE}_{\{B\}}$, $e$ is special in $G$.  If $G^{*}$ is not
$B$-free, then as $G$ is $B$-free, it follows that one of the other two
edges of $t$ is special in $G$.
\qed

\section
7.  New universal graphs.

We gave a general construction in Theorem 4 which produces finite sets
$\CC$ of connected constraints for which $T^{*}_{\CC}$ is
$\aleph_0$-categorical and hence, in particular, there is a universal
$\CC$-free graph.  We now generalize this.

\proclaim{Theorem 5}
Let $\CC$ be a finite set of finite connected graphs such that 
$T^{*}_{\CC}$ is $\aleph_0$-categorical.  Let $\cal H$ be a finite
set of finite connected graphs which is closed under homomorphic image.
Then $T^{*}_{{\CC}\cup {\cal H}}$ is $\aleph_0$-categorical.
\endproclaim

\proof:
Let $G\in {\EE}_{{\CC}\cup{\cal H}}$, and $A\subseteq G$ finite.  We
must bound $|\acl(A)|$ in terms of $|A|$ and apply Theorem 3.  
Let $G'\in {\EE}_{\CC},G\subseteq G'$.  Let 
$B=\acl_{G'}(A)$, a finite set of size bounded by a function of $|A|$.  It
suffices to show that $B_0=B\cap G$ is algebraically closed in $G$. 

Suppose the contrary, by Lemma 5 we have some $C\in {\CC}\cup {\cal H}$
and
a homomorphic image $C'\subseteq G$ so that $C$ embeds in the free amalgam
over $B_0$ of $|C|$ copies of $C'$.  If $C\in {\cal H}$ then 
$C'\in {\cal H}$, contradicting the assumption that $G$ omits 
${\CC}\cup {\cal H}$. Hence
$C\in {\CC}$.

Now we consider $B_0,B$, and $C'$ in $G'$.  As $G'$ omits $C$, $G'$ does
not contain the free amalgam of $|C|$ copies of $C'$ over $B_0$ or over
$B$ $(C'\cap B=C'\cap B_0)$.  Then by Park's theorem, $B$ is not
algebraically closed in $G'$, a contradiction.
\qed

It is of course trivial to produce examples of set of constraints to play
the role of $\cal H$ here, but we will want to consider a number of
concrete constructions, particularly with a view toward keeping 
$|{\cal H}|$ small.  This will require some preliminary observations.

\noindent{\bf Remarks}

1. We do not in fact require $\cal H$ to be {\sl closed} under the 
formation of homomorphic images. What is needed is the following: if 
$C'$ is a homomorphic image of 
$C\in {\cal H}$, then $C'$ contains an element of
$\cal H$.  In the future we will take this condition as the 
{\sl definition} of ``closure under homomorphism".  

2. In particular if $A$ is a finite connected graph we will write 
${\cal H}_0(A)$ for the set of all homomorphic images of 
$A$ and ${\cal H}(A)$
for the set of minimal elements of ${\cal H}_0(A)$ (with respect to
embeddings as subgraphs).  For example, if
$A$ is a cycle of odd length $2N+1$, then ${\cal H}(A)$ consists of odd
cycles of length $2n+1, n\leq N$. Similarly, if $A$ is a bipartite graph
containing at least one edge, then ${\cal H}(A)=\{K_{2}\}$.  More general,
for any finite connected graph $A$, ${\cal H}(A)$ contains a unique
complete graph $K_{n}$, with $n=\chi (A)$ the chromatic number.  Thus one
only gets new examples by considering graphs of chromatic number $\chi$
which do not contain the complete graphs $K_{\chi}$. In this case 
$|{\cal H}(A)|\geq 2$.

\proclaim{\bf Definition 7}
1.  Let $A_{1},A_{2}$ be two graphs.  Then $A_{1}\times A_{2}$ is the
graph with vertex set $V(A_{1})\cup V(A_{2})$, and whose edges are those
of $A_1$ and $A_2$ together with all pairs $(u,v)$, where $u\in A_{1},
v\in A_{2}$ or vice versa.

2.  Let ${\CC}_1$, ${\CC}_2$ be two sets of graphs.  Then 
${\CC}_{1}\times {\CC}_{2}=\{A_{1}\times A_{2}: 
A_{1}\in {\CC}_{1}, A_{2}\in {\CC}_{2}\}.$
\endproclaim

\noindent{\bf Remark}

${\cal H}_0(A_{1}\times A_{2})={\cal H}_0(A_{1})\times 
{\cal H}_0(A_{2})$.

\proclaim{\bf Example 11}
With $M,N$ fixed integers, the class of graphs omitting $C_{2m+1}\times
C_{2n+1}$ for $m\leq M, n\leq N$ has a universal graph; in particular for
$M=0$, this is the class constrained by forbidding ``wheels" 
$\{K_{1}\times C_{2n+1}:n\leq N\}$. 
\endproclaim

Another family of well-behaved examples is generated by application
of a construction used by Mycielski to generate triangle free graphs of
arbitrary high chromatic number, where Mycielski would begin with $K_2$,
we
substitute $K_n$, getting the following graphs, which we call $M_n$. Let 
$V(M_{n})= \{0\}\cup (\{1,2,\cdots,n\}\times \{0,1\})$ and 
set 
$A_n=\{1,2,\cdots,n\}\times \{0\}, B_{n}=\{1,2,\cdots,n\}\times \{1\}.$
Edges are defined as follows.  The vertex $0$ is adjacent to the vertices
of $A_n$ and no others; $M_n$ induces a complete graph on $B_n$, and no
edges on $A_n$; and the vertices $(i,0)$ and $(j,1)$ are adjacent if and
only if $i\not=j$.  This graph arises by applying Mycielski's construction
to $B_{n}$ (i.e. $K_n$).  To have a more suggestive notation we write
$a_i$ for $(i,0)$ and $b_i$ for $(i,1)$.

\proclaim{Lemma 7}
${\cal H}(M_{n})=\{K_k\times M_{n-k}:k\leq n, k\not=n-1\}.$ In
particular $\chi (M_{n})=n+1$ and $|{\cal H}(M_{n})|=n.$
\endproclaim

\proof:
Let ${\CC}=\{K_{k}\times M_{n-k}:k\leq n\}.$ Then 
${\CC}\subseteq {\cal H}_0(M_{n})$.  To see this, identify the
vertices $a_i$ and $b_i$ for $i>n-k$.  In particular for $k=n$ we have 
$K_{n}\times M_0=K_{n+1}\in {\CC}$, so 
$\chi (M_{n})\leq n+1$, and $\chi (M_{n})\geq n+1$ by Mycielski's argument
[BM, $\S 8.5$].

Any homomorphic image of $M_{n}$ other than those listed will involve
either the identification of the vertex $0$ with a vertex of $B_n$, or the
identification of vertices in $A_n$.  In either case the resulting
homomorphic image contains $K_{n+1}$ by inspection.  Thus the minimal
homomorphic images of $M_n$ belong to 
${\CC}: {\cal H}(M_{n})\subseteq {\CC}.$ Furthermore for 
$k=n-1, K_{k}\times M_{n-k}=K_{n-1}\times M_{1}\supseteq K_{n-1}\times 
K_{2}=K_{n+1}$, so
$K_{n-1}\times M_{1}\notin {\cal H}(M_{n})$.

It remains to be shown that the graphs $K_{k}\times M_{n-k}$ for 
$k\leq n, k\not= n-1$ are incomparable; this will complete the
characterization of ${\cal H}(M_{n})$.

Suppose therefore that $K_{k}\times M_{n-k}$ embeds in 
$K_{l}\times M_{n-l}$ with $0\leq k, l\leq n$ and 
$k,l\not= n-1, k\not= l$.  As 
$|K_{k}\times M_{n-k}|\leq |K_{l}\times M_{n-l}|$ we have $k\geq l$.
The case $k= n$ may be eliminated by inspection.  Accordingly we assume 
$0\leq l < k \leq  n-1$.   

Let $f:K_{k}\times M_{n-k}\rightarrow K_{l}\times M_{n-l}$ be an
embedding.  
As $k>l$, fix $u\in K_{k}$ so that $f(u)\notin K_{l}$.  
Now $u$ has $2n-k$ neighbors in $K_{k}\times M_{n-k}$.  
If $f(u)\notin B_{n-l}$, then $f(u)$ 
has at most $n$ neighbors in $K_{l}\times
M_{n-l}$, forcing $2n-k\leq n$, a contradiction.  
So $f(u)\in B_{n-l}$ and as $u$ is adjacent to every 
other vertex of $K_k\times M_{n-k}$, 
$f[K_{k}\times M_{n-k}]$ does not contain the vertex labelled $0$ in
$K_{l}\times M_{n-l}$.  
However the graph resulting from deletion of this
vertex has chromatic number $n$, while $K_{k}\times M_{n-k}$ has chromatic
number $n+1$, a contradiction.
\qed

Examples of constraint sets $\CC$ allowing a universal graph with
$|{\CC}|=1$ are very rare, and indeed few examples are known with any
sharp bound on $|{\CC}|$.  We will consider the possibilities in the
case ${\CC}={\cal H}(A)$.  Evidently, if we require $|{\CC}|=1$ we
will have ${\CC}=\{K_{n}\}$ for some complete graph, which is one of
the oldest examples.  We can on the other hand produce a number of
new examples with 
$|{\cal H}(A)|=2$.  We note first the simple example
${\cal H}(C_{s}\times K_{n})=\{C_{s}\times K_{n}, K_{n+3}\}$.  It seems
possible {\sl a priori} that these are the only such examples, and we
therefore will give an additional construction, showing at least that it
will not
be easy to classify the cases with $|{\cal H}(A)|=2$.

\proclaim{Construction}
Let $G=r\cdot K_{n}+K_{m}$ be the disjoint sum of $r$ complete graphs
$K_{n}$, and one more, $K_{m}$, with $n\geq m\geq 1$ and either $r\geq 2$,
or $m\geq n-1$.

We will write 
$G=A_{1}+\cdots +A_{r}+B$ with $A_{i}\simeq K_n$, $B\iso K_m$.  Let
$m_0=min (m,n-1)$ and let 
$\Sigma = \{S\subseteq V(G):
\hbox{ $|S\cap A_i|=n-1$ for $1\leq i\leq r$ and $|S\cap B|=m_0$}
\}$.

Let $G^{*}\supseteq G$ be defined as follows:
$V(G^{*})=V(G)\cup \{v_{S}:S\in \Sigma \}$. $G^{*}$ induces $G$ on
$V(G)$. The $v_{S}$ for $S\in \Sigma$ form an independent set, and the
neighbors of $v_{S}$ in $V(G)$ are the elements of $S$.
\endproclaim

\proclaim{\bf Example 12}
For the simplest example, take $r=1,n=2,m=1$.  Then $G=K_{2}+K_{1}$, and 
$G^{*}\simeq C_{5}$.
\endproclaim

\proclaim{Lemma 8}
$G^{*}$ defined above has chromatic number $n+1$.
\endproclaim

\proof:
One can color $G^{*}$ with $n+1$ colors by first coloring $G$ with $n$
colors, and using the last color for all remaining vertices.

On the other hand, if $G^{*}$ is colored with $n$ colors, then all $n$
colors occur in $A_{1}$.  Fix $b\in B$.  For each color $c$ fix $S(c)\in
\Sigma$ so that $b\in S(c)$ and $S(c)\cap A_{1}$ consists of those
vertices not of color $c$.  Thus $v_{S(c)}$ must have color $c$, so $b$
does not have color $c$.  Therefore $b$ cannot be colored. 
\qed

Now we give additional examples of constraint families $\CC$ such that
the algebraic closure is trivial in ${\EE}_{\CC}$ (i.e..
$\acl(A)=A$), and $|{\CC}|=2$.

\proclaim{Proposition 2}
For $r,m,n$ with $n\geq m\geq 1$ and either $r\geq 2$ or $m\geq n-1$, and
for $G^{*}$ defined as above,
${\cal H}(G^{*})=\{G^{*}, K_{n+1}\}$.  In particular
$|{\cal H}(G^{*})|=2$.
\endproclaim

\proof:
Evidently $G^{*}$ does not contain $K_{n+1}$.  It suffices now to prove
that any proper homomorphic image of $G^{*}$ does contain $K_{n+1}$.

Let $h:G^{*}\rightarrow H$ with $h(u)=v$, for some $u,v\in
V(G^{*}),u\not=v$.  Note that $u,v$ are not adjacent.  We consider cases.

Case 1. $u\in A_{i}$ for some $i$, and $v\in G$.

Take $S\in \Sigma$ with $S\cap A_{i}=A_{i}-\{u\}$, $v\in S$.  Then the
induced graph on $\{v_{S}\}\cup A_{i}$ is isomorphic to $K_{n+1}$ in $H$.

Case 2. $u\in G$, for some $i$, and $v\notin G$.

Let $A^*=A_i$ or $B$ be the component of $G$ containing $u$.
As $r\geq 2$ or $m\geq n-1$, we may choose 
$A\subseteq G, A\simeq K_{n-1}$,with $A\intersect A^*=\emptyset$,
and $v$
adjacent to all vertices of $A$.  Take $S$ with $u\in S$, $A\subseteq S$.
Then in $H$, the induced graph on $\{v\}\cup \{A\}\cup \{v_{S}\}$ is
isomorphic to $K_{n+1}$.

Case 3. $u,v\notin G$ 

Let $u=v_{S},v=v_{T}$. Let $A$ be a connected component of $G$ such that 
$A\cap S\not= A\cap T$.  Then in $H$, the induced graph on 
$\{u\}\cup A$ is isomorphic to $K_{n+1}$.
\qed

\section
8. Another Universal Graph.

The main result of this section is that for the
graph $C=T_1+_\cdot T_2+_\cdot P_n$ consisting of two triangles 
$T_{1},T_{2}$ with
exactly one common vertex and a path
of length $n$ starting from a non-common vertex in one of these triangles,
the theory $T_C^*$ is $\aleph_0$-categorical.
Here we use the {\sl ad hoc} notation $+.$ for an almost disjoint sum
with one (specific) pair of vertices identified.  (We write $+_v$ when the
common vertex $v$ needs to be specified).
This depends on an analysis of algebraic closure much of which is valid more
generally and may be useful in the analysis of other candidates for
membership in $\UU_0$ (defined in the introduction).

We will assume throughout that  $\CC$ is a finite set
of finite connected graphs.  Furthermore $G$ denotes an $\aleph_0$-saturated 
graph in ${\EE}_{\CC}$. We use the term ``weak embedding'' for the
ordinary graph theoretic
embedding (as opposed to a strict embedding, which is an isomorphism with an 
induced subgraph).

\proclaim{\bf Definition 8}
\item {1.} For $A\subseteq H\subseteq G$ with $H$ finite, we say that
$H$ is {\sl free} over
$A$ if there is an embedding of infinitely many copies of $H$ in $G$ over $A$, 
disjoint
over $A$ (this means that the intersection of any two copies  is $A$).

\item {2.} For $A\subseteq H\subseteq G$  with $H$ finite, 
$\cl(A;H)$ is the union of $A$ with all sets $B$ such that:
\itemitem{2.1} $H$ is free over $A\cup B$;
\itemitem{2.2} $B$ is minimal subject to $2.1$.

\item {3.} Let $\FF$ be a collection of pairs $(A,H)$ of finite graphs
with
$A\subseteq  H$. 
Then for any $X\subseteq G$, $\cl(X;{\FF})$ is
the union of all sets of the form $\cl(A;H)$ 
where $A\subseteq X$, $H\subseteq G$,
and $(A,H)$ is isomorphic to a pair in $\FF$.

\item{4.} With $\FF$ as in $3$, we say that $\FF$ is a {\sl base}
for $\acl$ if for all $X\subseteq G$
we have: $X=\cl(X;\FF)$ if and only if $X=\acl(X)$.

\item{5.} A graph $H$ is {\sl solid} if every induced 2-connected subgraph
of $H$ is complete. 
\endproclaim

We may now state the main results:

\proclaim{Proposition 3.}
For any pair $A,H$ of finite graphs with 
$A\subseteq H\subseteq G$, we have
$\cl(A; H)\subseteq \acl(A)$. 
Hence for any collection $\FF$ of 
pair $(A,H)$ of finite graphs with
$A\subseteq H$, and any $X\subseteq G$,
$\cl(X;{\FF})\subseteq \acl(X)$. 
\endproclaim

\proclaim{Proposition 4.}
If $\FF$ is a finite set of pairs $(A,H)$
of finite graphs with
$A\subseteq H$, and $X\subseteq G$ is finite, then 
$\cl(X;\FF)$ is finite.
\endproclaim

\proclaim{Proposition 5.}
Let ${\FF} =\{(A,H):A\subseteq H \subseteq G$
and for some $C\in \CC$, $H$ embeds weakly in $C$ as a proper subgraph
of $C$$\}$.
Then $\FF$ is a base for $\acl$.
\endproclaim

\proclaim{Proposition 6.}
If $\CC$ consists of solid graphs, and if 
$\FF$ is the collection of pairs
$(\{a\},H)$ for which $a\in H$, $H$ embeds properly in some $C\in \CC$,
and $H-\{a\}$ is a 
connected component of $C-\{a\}$, 
then $\FF$ is a basis for $\acl$. In particular:
$$\acl(X)=\bigcup_{a\in X}\acl(a)\leqno(*)$$
for $X\subseteq G$. 
\endproclaim

We do not know exactly when the ``unarity'' condition (*) holds; it might be 
useful to determine this.
If we take the union of a collection of solid graphs and a collection closed 
under homomorphic image, 
then the same property holds since $\acl$ is unchanged. However, 
if ${\CC}=\{C\}$ consists of a single forbidden subgraph, then $(*)$ is 
equivalent to the solidity
of $C$.

\proclaim{\bf Definition 9}
\endproclaim
The next statement requires a more delicate partial closure operation, for 
use with the particular
graph $C=T_1+_{u_2}T_2+_{y_0}P_n$ referred to above.
Let $\FF$ be the set of pairs $(\{a\},P)$
for which $P$ is a path of length at most $n$ with $a$ an endpoint.
For $X\subseteq G$ let $\cl^*_C(X)$ be the union of $\cl(X;\FF)$ with 
\item{1.} all sets of the form
$\cl(\{a\};H)$ for which: $a\in X$; 
$a$ lies in some copy of $T_1+_{\cdot}T_2$ with $a$ not the
common vertex of the two triangles; 
$H\simeq T+_{\cdot}P$, 
the free amalgam of a triangle $T$ 
with a path of length at most $n$; and
\item{2.} the set of all points $b$ lying in $\cl(\{a\};H)$ with $a\in X$, 
$a\in H\simeq T+_{\cdot}P$, and either $a,b$ belong to a triangle, or $b$
lies in some copy of $T_1+_{\cdot}T_2+P$ with $b$ not the
common vertex of the two triangles, and with $P$ a path. 

\proclaim{Proposition 7.}
Let $C=T_1+_{u_2}T_2+_{y_0}P_n$ 
be the graph referred to above, obtained
by amalgamating two triangles $T_1$, $T_2$, and a path $P_n$ of length $n$, 
over two distinct
points of $T_2$. Then for $X\subseteq G$ and $\cl_C^*$ as defined above,
if $X=cl_C^*(X)$ then $X=\acl(X)$.
\endproclaim

\proclaim{Proposition 8}
Let $C=T_1+_{u_2}T_2+_{y_0}P_n$ be the graph referred to 
in the previous Proposition. Then the theory $T_{\CC}^*$ is
$\aleph_0$-categorical,
and thus there is a universal $C$-free graph.
\endproclaim

\proof\  of Proposition 3:
We consider $A\cup B\subseteq H$ with $H$ free over $A\cup B$, and
with $B$ minimal subject to this condition
(in particular $A\cap B=\emptyset$). 
If $B\not\subseteq \acl(A)$, then there are infinitely many copies
$(B_i,H_i)$ of $(B,H)$ embedded 
as induced subgraphs of $G$, with the $B_i$ distinct and the $H_i$ free over $B_i$. Without loss of generality the $B_i$ form a 
$\Delta$-system  with common part $B_0$. As the $B_i$ are disjoint over $B_0$ 
and each $H_i$ is free over $A\cup B_i$, 
$H$ is free over $A\cup B_0$. This contradicts the minimality of $B$.
\qed

\proof\  of Proposition 4:
It is easy to see that $\cl(X;\FF)$  is a definable set, and as it is
contained in $\acl(X)$, and $G$ is $\aleph_0$-saturated, 
it is finite. For the definability it suffices to check the definability of 
``free over'';
but we can replace the requirement of infinitely many disjoint copies of $H$ 
by $k$ disjoint copies, where $k=\max\{|C|:C\in \CC\}$,
since $G\in {\EE}_{\CC}$. 
\qed

\proof\  of Proposition 5:
Let $\FF=$
$$\{(A,H):\hbox{ $A\subseteq H \subseteq G$
and for some $C\in \CC$, $H$ embeds weakly in $C$ as a proper subgraph
of $C$}\}$$
Let $X\subseteq G$, and assume that $X=\cl(X;\FF)$. We claim
$X= \acl(X)$. We 
may suppose that $X$ is finitely
generated. 

Suppose $X\ne \acl(X)$. Then as $G$ is $\aleph_0$-saturated, if we form 
$G(2)=G_1+_XG_2$ with $G_1$, $G_2$ isomorphic to
$G$ over $X$, then $G(2)\notin \GG_{\CC}$, and thus there is a
weak
embedding $h:C\hookrightarrow G(2)$ for some
$C\in \CC$. Let $H_2=h[C]\cap G_2$, and let $H_1$ be the image of
$H_2$ in $G_1$ under the given isomorphism.
As $X=\cl(X;\FF)$ and the pair $(X\cap h[C], H_1)$ lies in
$\FF$,
$H_1$ is free over $X\cap h[C]$, and hence
can be embedded in $G_1$ disjoint from $G_1\cap h[C]$ over 
$X\cap h[C]$. Defining $h'$ to agree with
$h$ off $h^{-1}[G_2]$ and with this new embedding on $h^{-1}[G_2]$, we have an 
embedding of $C$ into $G_1$. 
As $G_1\in \GG_{\CC}$ this is a contradiction.
\qed

The next two proofs will be somewhat similar to the foregoing, and very 
similar to one another.

\proof\  of Proposition 6:
$\CC$ consists of solid graphs and 
$\FF$ is the collection of pairs
$(\{a\};H)$ for which $a\in H$, $H$ embeds properly in some $C\in \CC$,
so that $H-\{a\}$ is a 
connected component of $C-\{a\}$. The proof that follows will allow us to 
replace $\FF$ by a slightly more
restricted family which will be defined below.

We take $X\subseteq G$ finitely generated (with respect to this closure
operation) and we suppose that
$X=\cl(X;\FF)$ but $X\ne \acl(X)$, so that after forming
$G(2)=G_1+_XG_2$ as in the previous argument,
we have an embedding $h:C\hookrightarrow G(2)$ for some $C\in\CC$.
We associate to $C$ the tree $T$ whose vertices 
correspond to the 2-connected components of $C$, with edges between components 
which either meet
or are connected by some edge of $C$. We will denote the vertices of $T$ by 
$t,t'$ and the like, and the component
of $C$ corresponding to a vertex $t$ of $T$ will be denoted $C_t$. Now pick an 
arbitrary vertex $0$ of $T$, 
and take it as a root for $T$. Now $T$ can be viewed as partially ordered set 
with minimum $0$. 
For $t\in T$ let $T^t=\{t':t'\ge t\}$ and let 
$C^t=\bigcup_{t'\in T^t} C_{t'}$. 

We will replace the set $\FF$ considered above by the subset of pairs 
$(\{a\};H)$ for which
for some $t>0$ either:
\item{1.} $a\in C_t$ and $H-\{a\}$ is a component of $C^t-\{a\}$; or
\item{2.} $t$ is a successor of a node $t^-$, $a\in C_{t^-}$, and 
$H-\{a\}$ is $C_t$.

\noindent Note that in the first case, typically $H=C^t$; this holds for 
example if $|C_t|>1$.

Of course, with this modification we still have $X=\cl(X;\FF)$.

Now let $t$ be maximal in $T$ such that $h[C^t]$ does not lie in either factor 
$G_1$ or $G_2$ of $G(2)$.
We may suppose that $h[C_t]\subseteq G_1$. 
It will suffice to replace $h:C\rightarrow G(2)$ on $C^t$ by $h'$
agreeing with $h$ on $(C-C^t)\cup C_t$, so
that $h'[C^t-C_t]\subseteq G_1$; repeating this operation eventually
produces an embedding of $C$ into
$G_1$, and a contradiction.  

We may break this down two steps further. First, for each successor node
$t'$ of $t$ for which $h[C^{t'}]\subseteq G_2$,
it suffices to find an embedding $h':C\rightarrow G(2)$ agreeing with $h$
on 
$C-C^{t'}$ and taking $C^{t'}$ into
$G_1$. For the second step, first choose a vertex $a\in X$ as follows: if 
$C_t{\cap C}_{t'}\neq \emptyset$,
let $a$ be the unique vertex common to both components. Otherwise, take a 
pair of vertices $u\in C_t$ and
$v\in C_{t'}$ with $u,v$ adjacent in $C$, and let $a=u$ if this is in $X$,
and $a=v$ otherwise. With these choices, $a\in X$. Now we adjust $h$ on
$\{a\}\cup C^{t'}$ by making separate adjustments on each subgraph
$H$ containing $a$ such that $H-\{a\}$ is a connected component of 
$C^{t'}-\{a\}$. 

At this point the pair $(\{a\};H)$ under consideration is one of the pairs 
which we have put in $\FF$.
As $X=\cl(X;\FF)$, we can embed $H$ freely into $G_1$ over $a$ and
arrive at 
the desired modification of $h$.
Iterating this construction over all such components and all such nodes $t'$, 
and then over all suitable $t$, we will
reach a contradiction.
\qed

\proof\  of Proposition 7:
We now deal with the particular case $C=T_1+_{u_2}T_2+_{y_0}P_n$, whose 
vertices we label as follows:
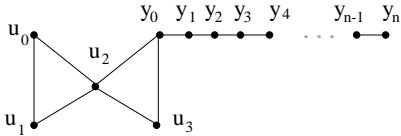
\epsfxsize= 15  em 
$$\epsfbox{689-pic.eps} $$

We follow exactly the same line as in the previous proof. Now the tree $T$ is 
a path of length $n+2$ whose
first node corresponds to the first triangle $T_1$; take this node as a root 
and use the corresponding set
$\FF$ of pairs:
$$\eqalign{
&(\{a\};P)\hbox{ with $a$ the initial vertex of a path $P$ of length at
most $n$}\cr
&(\{u_2\},T_2+_{y_0}P_n)\cr
}
$$
This is almost what we want, except that the second possibility is somewhat 
more generous than we wish to allow.
Accordingly, we will now consider the corresponding part of the previous 
argument more carefully.
This occurs when the vertex $t$ is the root and $t'$ corresponds to the 
triangle $T_2$, $a=u_2$, and
we wish to embed $H=T_2+_{y_0}P_n$ into $G_1$ over $a$ disjoint from the 
image of $T_1$.
That is, we have an embedding $h:C\rightarrow G(2)$ with $h[T_1]\subseteq
G_1$; 
$h[H]\subseteq G_2$, 
(so $h(u_2)\in X$) and
we assume toward a contradiction that any embedding of $H$ into $G_1$ meets 
$h[T_1]$. 
Let $b=h(u_3)$ and $c=h(y_0)$.

If $b,c\in X$ 
then it suffices to embed $P_n$ into $G_1$ correctly, and this we have 
already dealt with.
If neither $b$ nor $c$ lies in $\cl(\{a\};H)$ then 
$T_2$ can be embedded freely in $G_2$ over $a$, which produces a copy 
of $C$ in $G_2$ since we already have $T_2+_{y_0}P_n$ embedded in $G_2$.

Thus we are left with the cases in which $b$ or $c$ lies in $\cl(\{a\};H)$, 
and in particular lies in $X$,
and the other vertex is not in $X$.

Suppose $c\in \cl(\{a\};H)$ and  $b\notin X$. Let $B=\cl(\{a\};H)$. Then we 
can embed $B$ freely in $G_1$
over $X$, and then continue to embed $P_n$ freely in $G_1$ over $y_0$. This 
produces the desired embedding
of $C$ in $G_1$ (since ``freely'' means: without any undesirable 
identifications).

Finally, suppose $b\in \cl(\{a\};H)$ and $c\notin X$. In particular $c$ is not 
in $\cl(\{a\};H)$ and hence
there are infinitely many triangles containing $a,b$.
Let $B=\cl(\{a,b\};H)$. We will show that $B\subseteq X$. 
Take $u\in B-\{a,b\}$, and set $B'=B-\{u\}$. 
Let $G'$ be the free amalgam of $G_2$ with a large number of copies of 
$H$ over $B'$. Then $C$ embeds in $G'$. 

Suppose this embedding involves a triangle $T=\{a,b,c'\}$ lying in one of 
the copies of $H$. 
Then there is a triangle $T_0$ meeting $\{a,b\}$ in a single vertex. 
If this triangle contains $c$, then choosing $c^*\in G_2$ not in $T_0\cup 
h[P_n]$, so that
$\{a,b,c^*\}$ lie on a triangle, we get an embedding of $C$ into $G_2$. So
$T_0$ meets $\{a,b,c\}$ in one vertex.
If $T_0$ contains $u$ then by definition $u\in cl_C^*(X)=X$.
So suppose it does not contain $u$. As $T\cup T_0$ is part
of a copy of $C$ embedded in $G'$, this embedding also involves a path $P$ of 
length $n$ attached
to $T$ or $T_0$, and not at their common point. If the path is attached to 
$T_0$, then replacing $c$ by a point
$c^*$ for which $a,b,c^*$ forms a triangle and $c^*$ lies off $T_0\cup P$, 
again $C$ embeds in $G_2$, a contradiction.
So $P$ is attached to $T$. $P$ is broken into various connected components 
by its intersection with $B'$. We will
alter the embedding so that $P$ becomes a path of length $n$ attached to $c$ 
and otherwise disjoint from
$T_0\cup \{a,b,c\}$. Those segments which lie in $G_2$ may be left as they 
are. The remainder lie
in copies of $H$, are attached at one or two points of $B'\cup \{c'\}$,
and correspond to segments in $G_2$ 
which are either free over $B'$, or contain the point $u$. As $u$ does not lie 
on $T_0$, if a segment
corresponding to one containing $u$ occurs, it may be replaced by two segments 
in $G_2$ joined at $u$, 
and free over $B$. Thus by choosing the embedding of $P$ carefully, one may 
embed $C$ in $G_2$, a contradiction.

Therefore in our original embedding of $C$ into $G'$, the copy of 
$T_1+_{\cdot}T_2$ embeds in $G_2$ and part
of the path $P$ is embedded in various copies of $H$ amalgamated over $B'$.
Again we can alter most of the embedding of $P$ to go into $G_2$, apart from 
segments which
correspond to segments in $G_2$ lying between two successive
points of $B'$, with the vertex $u$ on the segment. If any such segment 
actually occurs, it means that
in $G_2$, $u$ lies on some graph of the form $T_1+_{\cdot}T_2+_{\cdot}P$ with 
$P$ a path. Thus again 
$u\in cl^*_C(\{a\})$. 
\qed

\medskip
For the next proof we will require an auxiliary result which will be seen to 
contain useful
information about algebraic closure in the case at hand.
Let $C=T_1+_{u_2}T_2+_{y_0}P_n$ 
be the graph referred to in Proposition 8.
In particular, $n$ is fixed.

\proclaim{Lemma 9}  
Let $G$ be a graph, $u$ a vertex of $G$, and suppose
that there are two disjoint paths of length $5n$ originating at $u$, as well 
as an embedding of some 
subgraph $H$ of $C$ of the form $C=T_1+_{u_2}T_2+_{y_0}P_k$  with 
$0\le k\le n$, embedded
with $u$ as the terminal vertex of $P_k$. Then $C$ embeds in $G$.
\endproclaim

\proof:
Let $v$ be the vertex in $G$ corresponding to the vertex $u_2$ of $H$. Let 
$P$ be one of the two given
paths, which does not contain $v$. Then $P$ is broken into at most $5$ 
connected components by its intersection
with the vertices of $T_1+_{u_2}T_2$ 
(as embedded in $G$), and one of these components has length at least $n$.
Thus, if this intersection is nonempty, then $C$ embeds in $G$. 

Suppose $P$ is disjoint from the image of 
$T_1+_{u_2}T_2$ in $G$. Let $y$ be the first vertex of the path $P_k$ 
(starting from the vertex $y_0$ in $T_2$)
which corresponds under the embedding to a vertex of $P$. Then on removal
of $y$ from $P$, one of the components
has length at least $n$, and hence we again have an embedding of $C$ into $G$.
\qed

\proof\   of Proposition 8:
We now wish to show that for $X$ finite, $\acl(X)$ is finite. We define 
inductively:
$X_0=X$, $X_{i+1}=cl_C^*(X)$, and we need to show that this process 
terminates. 
Suppose in fact that it goes on for $k$ stages with $k$ substantially larger 
than $10n$.
Define a sequence of points $a_i\in X_{i}\setminus X_{i-1}$ for  $i<k$ 
by downward induction 
so that $a_i\in cl_X^*(a_{i-1})$ for all $i$. 
The point $a_{k-1}$ is selected arbitrarily, and given $a_i$,
as it lies in $X_i$ it lies in $\cl_C^*(a_{i-1})$ for some 
$a_{i-1}\in X_{i-1}$, and this element lies
outside $X_{i-2}$ since $\cl_C^*(a_{i-1})$ is not contained in $X_{i-1}$. 

We claim: 
$$\hbox{ The elements $a_i$ ($i<k$) can be selected so that they lie on a 
path of length at least $k-1$.}\leqno(*)$$ 
Again, proceed by downward induction, building up a finite path 
$Q_i$ with endpoint $a_i$
as we go along in such a way that $Q_i\cap X_i=\{a_i\}$.

Suppose $b=a_i$ has been chosen and pick some $a\in X_{i-1}$ so that 
$b\not\in cl_C^*(a)$.  
Suppose first that $b\in \cl(\{a\};P)$ for some path $P$, with
$a$ an endpoint of $P$.  Let $B=\cl(\{a\};P)$ and let $a'$ be the vertex of
$B\cap X_{i-1}$ on the segment from $a$ to $b$ which is closest to $b$;
possibly $a'=a$.  Let $P'$ be the subpath of $P$ with initial vertex $a'$,
passing through $b$, and $B'=B\cap P'$.  We know that $P'$ is free over
$B'$ and easily $B'$ is minimal with this property. Thus $B'\subseteq
\cl(\{a\};P')$ and we may take $a_{i-1}=a'$. The remaining elements on the
segment $(a',b)\cap B'$ are outside $X_i$ and as $[a',b]$ is free over its
intersection with $B'$, we may attach to $Q_i$ a path from $a_{i-1}$ to
$a_i$ which meets $X_{i-1}$ only in $a_{i-1}$, and meets $Q_i$ only in
$a_i$; this produces the desired path $Q_{i-1}$. 

Now suppose that $b\in \cl(\{a\};H)$ with $H\simeq T+_{\cdot}P_n$, 
$a\in T$, $a\notin P_n$.
We can proceed in more or less the same way. If $b$ is a vertex of $T$ we can 
just take $a_{i-1}=a$ and
adjoin the edge $(a,b)$ to the path $Q_i$. Otherwise $b$ lies on the path 
$P_n$.  We consider 
$P=[a,b]$, the shortest path from $a$ to $b$ in $h$; this meets $T$ in two 
points.
Let $B=\cl(\{a\};H)\cap P$ and let $a'$ 
be the vertex of $B\cap X_{i-1}$ which is closest to $b$; possibly $a'=a$.
Let $P'$ be the subpath of $P$ with initial vertex $a'$, passing through $b$. 
If $a'\ne a$ we claim that $b\in \cl(\{a'\};P')$. This is seen as in the 
previous case.
Furthermore $P'$ is free over $B\cap P'$ so we may connect $a'$ to $b$ by
a path meeting $X_{i-1}$ in 
$a$ alone, and meeting $Q_i$ in $b$ alone. 

Thus we have $(*)$, and in particular if $b_i=a_{5n+i}$ and $i$ is not too 
large, we have 
two disjoint paths from $b_i$ of length $5n$ contained in $G$.
We claim that in this case every path of length $n$ originating at $b_i$ is 
free over $b_i$.
If this fails, then as $G\in {\EE}_C$ we must have a subgraph of $C$
consisting 
of the two triangles
and some initial segment $I$ (possibly of length 0) of the path $P_n$, 
embedded with $b_i$ as the terminal
point of $I$. This violates the previous lemma.
It follows that for $a=b_i$ and $b=b_{i+1}$ the relationship $b\in cl_C^*(a)$
is realized in the following way:
$b\in \cl(\{a\};H)$ with $H\simeq T+_{\cdot}P$ (a triangle amalgamated with
a path of length $n$), 
and either $a$ lies in some copy of $T_1+_{\cdot}T_2$ with $a$ not the 
central vertex, contradicting the
previous lemma, or $b$ lies in some 
copy of $T_1+_{\cdot}T_2+P$, with $P$ a path and  with $b$ not the central 
vertex, again contradicting the
previous lemma if $i+1$ is not too large, or finally: $a$ and $b$ lie on a 
common triangle.

Thus we may assume that we have 4 consecutive points $b_i$ ($i=1,2,3,4$) such 
that for $i=1,2,3$ the pair
$(b_i,b_{i+1})$ lies on a triangle $(b_i,b_{i+1},c_i)$. This gives an 
embedding of $T_1+_{\cdot}T_2$
into $\{b_i,c_j:1\le i \le 4, 1\le j\le 3\}$: if $c_1=c_3$ use the first and 
third triangles, while
otherwise these two triangles are disjoint and hence the second one meets at 
least one of them in a single
vertex.  This again violates the previous lemma since at least one of the 
$b_i$ occurs as a noncentral point
in the embedded copy of $T_1+_{\cdot}T_2$.

This contradiction completes the proof.

\section 9. Paths.

The existence of a universal $P$-free graph, when $P=P_k$ is a finite
path of length $k$ (and thus of order $k+1$) is established in [KMP].
The analysis given there yields good structural information and allows
further generalization, for example to categories of vertex colored
graphs, which will be of further use even in the case of graphs.
However it does not give realistic control over the sizes of
algebraic closures. Writing $c_k(n)$ for $\max|\acl(A)|$, where $A$
varies over sets of $n$ vertices in graphs $G$ belonging to
$\EE_{P_k}$, we would get an estimate of $c_k(n)$ of the form a tower
of exponentials of height about $k$, using [KMP]. However 
the analysis of \S8 yields:
$$c_k(n)=c_k(1)\cdot n$$
which is already fairly good, leaving open only the question of the
growth rate of $c_k(1)$ as a function of $k$, which turns out to be
an intriguing question. Consideration of circuits of length $k$, or,
for that matter, any hamiltonian graphs of order $k$, yields:
$$c_k(1)\ge k$$
and for low values of $k$ one may check $c_k(1)=k$. In fact the
following is open:

\proclaim{Problem. Is $c_k(1)$ equal to $k$ for all $k$?}
\endproclaim

In the remainder of this section we will prove:

\proclaim{Proposition 9} $c_k(1)< k^{3k^2}$.
\endproclaim

That is, we reduce a tower of exponentials to a single exponential, 
but fall far short of the linear bound which may hold.
This result requires a closer and more concrete analysis of the
operation of algebraic closure, which begins by simply following
through on the analysis given in \S8 more generally.

On the basis of Proposition 6 we can describe the algebraic closure
operation in $\EE_{P_k}$ as follows. Let $\FF$ be the collection of
pairs $(\{a\},Q)$ where $Q$ is a path of length at most $k-1$ and $a$
is an endpoint of $Q$. For $G\in \EE_{P_k}$ and $v\in V(G)$, 
define inductively: 
$$\cl_0(v)=\{v\}; \quad\cl_{n+1}(v)=\cl(\cl_n(v);\FF)-\cl_n(v)$$
if $\Union_{i=0}^{k-1}\cl_i(v)\ne \acl(v)$, one produces a
contradiction
by constructing a path of length $k$ by downward induction, beginning
with some $u\in \cl_k(v)$. Thus if we have a uniform estimate of the
form $|\cl_1(v)|\le N$ holding in $\EE_{P_k}$ (all $v$), then correspondingly
$$\sup |\acl(v)|\le \sum_{i\le k-1} N^i<N^k, 
\hbox{ (the supremum is over all $v\in G\in \EE_{P_k}$)}
$$

We will get such an estimate with $N=k^{3k}$.

\proclaim{Definition 10}
Let $G\in \EE_{P_k}$. 
\item{1.} For $u,v\in V(G)$, set $\omega(u,v)=$
$$\eqalign{
\sup\{m:&\hbox{ There are infinitely many paths of length $m$
connecting $u$ and $v$ in $G$,}\cr
&\hbox{ disjoint except for their endpoints}\}\cr
}
$$
This supremum is taken to be 0 if there is no such $m$. However
when $u$ and $v$ are adjacent the condition is considered to hold
with $m=1$, in a degenerate form.
\item{2.} For $u\in V(G)$, set $\omega(u,\infty)=$
$$\eqalign{
\sup\{m:&\hbox{ There are infinitely many paths of length $m$
	 in $G$ with $u$ as an endpoint,}\cr
&\hbox{ disjoint apart from $u$}\}\cr
}$$
\endproclaim

\proclaim{Lemma 10.}
Let $G\in \EE_{P_k}$, $v\in V(G)$, and suppose that
$Q$ is a path in $G$ originating at $v$, while $B\includedin Q-\{v\}$
is minimal such that $Q$ is free over $B\union \{v\}$. Write $B\union
\{v\}$
as a sequence $(v_0,v_1,\dots, v_l)$ in order along $Q$, beginning
with $v_0=v$. Then $\omega(v_i,v_{i+1})\ge 1$ for all $i<l$ and:
\item{1.} 
$\omega(v_i,v_j)<\sum_{i\le r<j}\omega(v_r,v_{r+1})$ for
$i\le j-2$; and
\item{2.} $\omega(v_i,\infty)<\sum_{i\le r<l}\omega(v_r,v_{r+1})+
\omega(v_l,\infty)$ for
$i<l$.
\endproclaim

\proof:
Since $Q$ is free over $B\union \{v\}$, we have
$\omega(v_i,v_{i+1})\ge 1$.

Condition (1) follows easily from the assumption that $Q$ is not free
over $(B\union \{v\})\setminus\{v_r:i<r<j\}$ for $i\le j-2$, and
condition
(2) follows from the assumption that $Q$ is not free over $(B\union
\{v\})\setminus \{v_r:r>i\}$. 
\qed

\proclaim {Definition 11}
Let $G\in \EE_{P_k}$ and let
$\vbar=(v_0,\dots,v_l)$ be a sequence of vertices of $G$.
\item{1.} $\vbar$ is a {\rm chain} if $l\ge 2$,
$\omega(v_i,v_{i+1})\ge 1$ for all $i$, and 
$\omega(v_i,v_j)<\sum_{i\le r<j}\omega(v_r,v_{r+1})$ for
$i\le j-2$.
\item{2.} Similarly, the formal sequence $(v_0,\dots,v_l,\infty)$ is called
an {\rm open chain} if $l\ge 1$ (so the length is at least 2) and
if it satisfies the same formal conditions, using $\omega(v_i,\infty)$
where called for. 
\item{3.} The {\rm virtual length} of a chain $\vbar$ is
$\sum_{r}\omega(v_r,v_{r+1})$, and the virtual length of an open chain
is defined similarly. We write $\lambda(\vbar)$ for the virtual length
of $\vbar$.
\endproclaim

In the proof of the next lemma we will need a result of Erd\"os and Gallai:

\proclaim{Fact [EG, Theorem 2.6]}
let $H$ be a graph with $n$ vertices and $e$ edges, in which there is no
path containing $l$ edges $(l\ge 1)$.  Then $e\leq {n(l-1)\over 2}$.
\endproclaim

\proclaim{Lemma 11.}
Let $G\in \EE_{P_k}$, $v\in V(G)$, and $A\includedin \cl_1(v)$.
Then there is a set $A'\including A$, with $|A'-A|<(k^3/4)|A|$, 
such that for any chain or
open chain $\vbar$ whose endpoints lie in $A\union \{\infty\}$, 
if $\vbar$ is not contained in $A\union \{\infty\}$ then it meets
$A'\setminus A$.
\endproclaim

\proof:
For each $a\in A$ choose one path originating at $a$, of maximum
length, and let $\AA_1$ be the set of paths chosen.  Let $\AA_2$ be a
maximal collection of chains whose endpoints lie in $A$, and which
are otherwise disjoint both from each other and from the paths in
$\AA_1$. We will take $A'=\Union(\AA_1\union \AA_2)$. There are a
number
of points to be verified. We will begin by verifying that $A'$ has the
desired property, then estimate its size.

Consider first a chain with endpoints in $A$, not wholly contained
in $A$. We may suppose then that only its endpoints lie in $A$.
By the choice of
$\AA_2$,
if this chain does not meet any path in $\AA_1$ in one of its interior
points, then it meets one of the chains in $\AA_2$. 

Now consider an open chain $\vbar$ 
originating at a vertex $a$ of $A$, and not
wholly
contained in $A\union \{\infty\}$. Then we may suppose that it meets
$A$ only at $a$, as otherwise we would replace it either by a shorter
open chain, or by a chain with endpoints in $A$. 
Let $L=\lambda(\vbar)$. Then $\omega(a,\infty)<L$, so there is a path
(an ``obstruction'') of length at least $k-L$ with $a$ as an endpoint.
Therefore there is such a path in $\AA_1$, and it is easy to see that
$\vbar$ meets that path at an interior point, as otherwise
one constructs a path of length $k$ in the ambient graph.

For cardinality estimates it will be convenient to take $k\ge 3$, as
we may.
We have $|\AA_1|\le |A|$ and to complete the analysis we will show:
$$|\AA_2|\le {(k-1)(k-1)\over 4}|A|,$$
from which our claim follows easily.

To make this estimate, we will estimate separately the number of
chains
in $\AA_2$ connecting two specified vertices, and the number of pairs
of vertices having such a connection. 

We begin with the latter point.
Consider the graph $\Gamma$ whose vertex set is $A$, and with edges
between pairs of vertices joined by one of the chains in $\AA_2$. 
Our claim is:
$$e(\Gamma)\le {k-1\over 4}|A|\leqno(*)$$
We claim that 
 $\Gamma$ contains no path of length $\ceiling{k/2}$;
this property then implies $(*)$ by the
result of Erd\"os and Gallai [EG] which was quoted above.
If $\Gamma$ had a path of length $\ceiling{k/2}$, 
that is to say a sequence of at least $k/2$
chains which are disjoint except at their endpoints, then these chains
fit together to form a path of length at least $k$ (extending the chains
``freely'' to their virtual lengths). Note that each chain has virtual length
at least 2 by definition, and the chains may be extended so that the
added vertices are distinct from each other and any vertices previously
considered.

Now we consider the number $\mu(a,b)$ of chains in $\AA_2$ which
connect
two fixed vertices $a,b\in A$. We claim $\mu(a,b)\le k-1$.  
To see this, let $l=\omega(a,b)$. 
Then for each chain $\vbar\in \AA_2$
which connects $a$ and $b$, we have $\lambda(\vbar)\ge l+1$.
Let $G_1$ be the free amalgam of the ambient graph $G$ with 
infinitely many additional paths of length $l+1$ connecting
$a$ to $b$. If $G_1\in \GG_{P_k}$ then as $G\in \EE_{P_k}$
and $\omega(a,b)\ge l+1$ in $G_1$, we find $\omega(a,b)\ge l+1$
in $G$, a contradiction.

Thus $G_1$ contains 
a path $P$ of length $k$. $P\setminus\{a,b\}$ consists of at most
3 segments, each of which lies either wholly in $G$ or in one of
the additional paths of length $l+1$ adjoined to form $G_1$.
Assuming $\mu(a,b)>0$, 
there is at least one path in $G$ of length $l+1$ joining $a$ and $b$.
Therefore
we may suppose that $P$ is chosen so that $P'=\P\intersect G$ contains
at least one of the segments of $P\setminus \{a,b\}$.
In particular $P\setminus P'$ consists of at most two segments.
We now count separately the chains $\vbar\in \AA_2$ connecting $a$ and $b$
which meet $P'\setminus\{a,b\}$, and those which do not. $P'$ contains
at most $k-2$ vertices and hence meets at most $k-2$ of the chains in
$\AA_2$. Furthermore $P'$ is disjoint from at most one chain in
$\AA_2$ which links $a$ and $b$, as two such chains could be extended
freely to give two disjoint paths of length $l+1$ joining $a$ and $b$,
into which the segments of $P\setminus P'$ could be copied, thereby
embedding
$P_k$ in $G$.
Thus there are at most $1+(k-2)=k-1$ chains in $\AA_2$ linking $a$ and $b$.
This completes our estimate.
\qed

\proclaim{Corollary}
For $G\in \EE_{P_k}$, $v\in G$, we have $|\cl_1(v)|<k^{3k}$.
\endproclaim
\proof:
Let $A_0=\{v\}$ and define inductively $A_{i+1}=A_i'$ in the sense of
Lemma 11 (this is not canonical, of course). In other words, choose
$A_{i+1}$ satisfying:
$$A_i\includedin A_{i+1},\quad |A_{i+1}\setminus A_i|<{k^3\over 4}|A_i|,$$
so that any chain or open chain $\vbar$ whose endpoints lie in $A_i\union 
\{\infty\}$ is either contained in $A_i\union \{\infty\}$, or meets 
$A_{i+1}\setminus A_i$.

Then $|A_{i+1}|\le k^3|A_i|$ and hence $|A_i|\le k^{3i}$ for all $i$.
On the other hand $\cl_1(v)\includedin A_k$ since each 
open chain $\vbar$ originating at $v$ will meet $A_{i+1}\setminus A_i$,
as long as it is not contained in $A_i$; and $\vbar$ has at most $k$
vertices.
\qed

\section 10. More examples of universal graphs.

We give two more examples of constraints allowing universal graphs.  
These are less complex than the family  treated in \S 8, and may allow
some further elaboration.

Consider first the constraint $C=K_n+_\cdot P_k$ consisting of a
complete graph with an attached path. It can be shown using either the
methods
of \S8 or those of [KMP] that $\EE_{C}$ is $\aleph_0$-categorical. 
Every connected component $G_0$ of a graph $G\in \EE_{C}$ either omits
$K_n$ and belongs to $\EE_{K_n}$, or contains a copy $K$ of $K_n$, in
which case the connected components of $G_0\setminus K$ omit $P_{2k}$,
and a structure theory for these can be given in the spirit of [KMP].
Alternatively, following the argument of \S8, one finds that a vertex
lying on a sufficiently long path has trivial algebraic closure.

Our second example is a slight generalization of the bow-tie, namely
$C=K_n+_\cdot K_3$, a complete graph attached to a triangle. A
detailed analysis of the algebraic closure operator in this
case will yield:
$$|\acl(A)|\le (n+1)|A|$$
for $A\includedin G\in \EE_{C}$. We will now give the details.

\proclaim{Definition 12}
Let $G\in \EE_C$.
\item{1.} For $A\includedin V(G)$, $A$ is {\rm special} if $G$
induces a complete graph on $A$, and one of the following occurs:
\itemitem{a.} $|A|=n$ and there is no $B\includedin V(G)$ of order $n$
such that $B\ne A$, $B\intersect A\ne \emptyset$, and $G$ induces
a complete graph on $B$; or
\itemitem{b.} $|A|=n-1$ and there are at least two vertices of $G$
adjacent to all vertices of $A$.
\item{3.} For $a\in V(G)$ set $a^*=$
$$\Union\{B:\hbox{$B$ is special, and
the graph induced on $\{a\}\union B$ is complete}\}$$
\endproclaim

\proclaim {Remark} \endproclaim
Let $G\in \EE_C$, $A\includedin V(G)$ of order $n$, and suppose that
$G$ induces a complete graph on $A$. Then $A$ contains a special
subset.

Our objective is to show that for $G\in \EE_C$ and $a\in V(G)$, we
have
$\acl(a)=\{a\}\union a^*$, and $|\{a\}\union a^*|\le n+1$, so that
by Proposition 6 of \S8 we may conclude $|\acl(A)|\le (n+1)|A|$ for all
$A$,
and in particular this constraint allows a (canonical) universal
countable graph.

\proclaim{Lemma 12} Let $G\in \EE_C$. 
\item{1.} If $A,B\includedin V(G)$ are of order $n$, and 
$G$ induces a complete graph on 
each, then either $A=B$, or $A\intersect B=\emptyset$, or
$|A\intersect B|=n-1$.
\item{2.} If $A,B\includedin V(G)$ are special and $A\intersect B\ne
\emptyset$, then either $A=B$ or $A\union B$ is contained in a set of
vertices of order $n+1$ on which $G$ induces a complete graph.
\endproclaim
\proof:
(1) holds by inspection.

For (2), if $|A|=n$ the claim holds by definition, so suppose
that $|A|=|B|=n-1$ and $A\ne B$. Let $A^+=A\union \{a\}$ and
$B^+=B\union \{b\}$ be two sets of vertices of order $n$ on which
$G$ induces complete graphs. We may suppose that they are chosen so
that $a\notin B$ and $b\notin A$. This forces $a=b$ as $|A^+\intersect
B^+|=n-1$. Now let $B^*=B\union \{b'\}$ be another choice for $B^+$,
so $b'\ne b$. If $b'\notin A$ then $G$ induces $C$ on $(B\union
\{b'\})\union (\{a\}\union A\setminus B)$, a contradiction. So 
$A\setminus B=\{b'\}$, which means that $G$ induces $K_n$ on $A\union
B$,
and induces $K_{n+1}$ on $\{a\}\union A\union B$. 
\qed

\proclaim{Corollary}
For $G\in \EE_C$, and $a\in V(G)$, we have $|\{a\}\union a^*|\le n+1$.
\endproclaim

\proclaim{Lemma 13} 
For $G\in \EE_C$, and $A\includedin V(G)$, the following are
equivalent:
\item{1.} $A=\acl (A)$;
\item{2.} $a^*\includedin A$ for $a\in A$.
\endproclaim
\proof:
($1\implies 2$):
Assume $A=\acl(A)$, $a\in A$,  and $B\includedin a^*$ 
is special, with $G$ inducing a complete graph on $\{a\}\union B$. 	
If $|B|=n$ this implies that 
$a\in B$ and $B$ is the unique
such set containing $a$, hence belongs to $\acl(a)\includedin A$. 
If on the other hand $|B|=n-1$ then either $B$ is unique, 
or else $B\union \{a\}$ is contained in a unique
complete
subgraph of $G$ of order $n+1$, by part (2) of the preceding lemma.
In either case $B\includedin \acl(a)\includedin A$.

($2\implies 1$):
We suppose $a^*\includedin A$ for $a\in A$, but $A$ is not
algebraically closed, and hence $C$ embeds into the free amalgam
$G_1+_A G_2$ of two copies of $G$ amalgamated over $A$. 
Let $C^*$ be the image of $C$. As $C^*\setminus A$ is disconnected,
there is a point $a\in C^*\intersect A$ which lies on a complete graph
$K$ of order $n$ and a triangle $T$ intersecting at $a$; we may take
$K$ to lie in $G_1$, and $T$ to lie in $G_2$. 

In particular $V(K)$ contains a special set $B$, and then
$B\includedin  a^*\includedin A$. As $G_2$ is $C$-free, this forces
$|\{a\}\union B|=n-1$, that is: $|B|=n-1$ and $a\in B$. 
Now let $T_1$ be the triangle in $G_1$ which corresponds to $T$ in
$G_2$ (via some isomorphism over $A$). Then $V(T)\intersect B=\{a\}$
but $T_1$ must have another vertex in common with $K$, as $G_1$ is 
$C$-free, and thus $T_1\includes V(K)\setminus B$ (which consists of a
single vertex). 
At the same time, as $B$ is special, there is another
complete
graph $K'$ of order $n$ in $G_1$ which contains $B$, and by the same
token the triangle $T_1$ contains $V(K')\setminus B$; so as $T_1$ 
is complete, it follows that the graph induced on $V(K)\union V(K')$
is complete of order $n+1$. But then $V(K)\includedin a^*\includedin A$
as all its subsets of order $n-1$ are special, and again
$C^*\includedin G_2$, a contradiction.
\qed

Thus as indicated above, we find $|\acl(A)|\le (n+1)|A|$, so
$\EE_C$ is $\aleph_0$-categorical and there is a universal $C$-free
graph. The case of a bouquet of two complete graphs, each of order
at least 4, has not been investigated and may well succumb to a
similar analysis. 
In any case, we believe that it should now be clear that the
classification of the class 
$\UU_0$, described in the introduction, is within reach, albeit 
this would involve some rather substantial computations in positive 
cases and 
some additional concrete constructions to cover the negative cases.
We emphasize that while the details would no doubt be tedious, the
result would be a reasonably well-founded conjecture as to the 
general solution  
of the problem of the existence of a countable universal graph,
for the class of graphs specified by prescribing {\sl any} (single)
finite connected forbidden subgraph.

\bigskip

\def\sameauthor{\vrule height.4pt depth 0pt width .4 true in}

\def\refindent{9 ex }
\def\ref[#1] {\hangindent \refindent\hangafter 1 
   \medskip\noindent \hbox to \refindent{[#1]\hfill}\ignorespaces}


\noindent{\bf References}

\ref[Ba] P.~Bacsich, {\sl The strong amalgamation property}, Colloquium
Mathematics, {\bf 33} (1975), 13-23.

\ref[BM] J.~A.~Bondy and U.~S.~R.~Murty, {\bf Graph Theory with Applications},
North-Holland, 1979.

\ref[CK] C.~C.~Chang and H.~J.~Keisler, {\sl \bf Model Theory}, North-Holland,
3nd ed., 1990.

\ref[ChK] G.~Cherlin and P.~Komj\'{a}th, {\sl There is no universal countable 
pentagon-free graph}, J.~Graph Theory, {\bf 18} (1994), 337-341.

\ref[CS1] G.~Cherlin and N.~Shi, {\sl Graphs omitting sums of complete 
graphs}, J.~Graph Theory, {\bf 24} (1997), 237-247.

\ref[CS2] \sameauthor, {\sl Graphs omitting a finite set of 
cycles}, J.~Graph Theory, {\bf 21} (1996), 351-355.

\ref[CST] G.~Cherlin, N.~Shi, and L.~Tallgren, {\sl Graphs omitting a bushy
tree}, J.~Graph Theory, {\bf 26} (1997), 203-210.

\ref[EG] P.~Erd\"os and T.~Gallai, {\sl On maximal paths and circuits of
graphs}, Acta Math., Acad.~Sci.~Hung.~{\bf 10} (1959), 337-356.

\ref[FK1] Z.~F\"{u}redi and P.~Komj\'{a}th, {\sl On the existence of countable 
universal graphs}, J.~Graph Theory, {\bf 25} (1997), 53-58.

\ref[FK2] \sameauthor, {\sl Nonexistence of universal graphs without some
trees}, Combinatorica, {\bf 17} (1997), 163-171.

\ref[GK] M.~Goldstern and M.~Kojman, {\sl Universal arrow free graphs}, 
Acta.~Math.~Hungary {\bf 73} (1996), 319-326.

\ref[HW] J.~Hirschfeld and W.~H.~Wheeler, {\sl \bf Forcing, Arithmetic,
Division Rings}, Springer-Verlag, 1975.

\ref[Ko] P.~Komjath, {\sl Some remarks on universal graphs}, Discrete
Math., to appear.

\ref[KMP] P.~Komj\'{a}th, A.~Mekler and J.~Path, {\sl Some universal graphs}, 
Israel J.~Math.~{\bf 64} (1988), 158-168.

\ref[KP1] P.~Komj\'{a}th and J.~Pach, {\sl Universal graphs without large
bipartite subgraphs}, Mathematika, {\bf 31} (1984), 282-290.

\ref[KP2] \sameauthor, {\sl Universal elements and the complexity of certain
classes of infinite graphs}, Discrete Math.~{\bf 95} (1991), 255-270.

\ref[Pa] J.~Pach, {\sl A problem of Ulam on planar graphs}, Eur.~J.~Comb.~
{\bf 2} (1981), 357-361.

\ref[Ra] R.~Rado, {\sl Universal graphs and universal functions}, Acta Arith.~
{\bf 9} (1964), 331-340.

\end